\documentclass[12pt]{amsart}

\usepackage{amscd,amsthm,amsfonts,amssymb,amsmath}
\usepackage{mathrsfs}
\usepackage[colorlinks=true]{hyperref}
\usepackage{fullpage}
\usepackage[all]{xy}
\usepackage{xcolor}

    \newcommand{\BA}{{\mathbb {A}}} \newcommand{\BB}{{\mathbb {B}}}
    \newcommand{\BC}{{\mathbb {C}}} 
    \newcommand{\BE}{{\mathbb {E}}} \newcommand{\BF}{{\mathbb {F}}}

     \newcommand{\BL}{{\mathbb {L}}}
     
     \newcommand{\BP}{{\mathbb {P}}}
    \newcommand{\BQ}{{\mathbb {Q}}} \newcommand{\BR}{{\mathbb {R}}}

     \newcommand{\BZ}{{\mathbb {Z}}}

     \newcommand{\CB}{{\mathcal {B}}}

     \newcommand{\CH}{{\mathcal {H}}}

    \newcommand{\CO}{{\mathcal {O}}} 
     \newcommand{\CR}{{\mathcal {R}}}

     \newcommand{\RN}{{\mathrm {N}}}

     \newcommand{\fp}{{\mathfrak{p}}}

     \newcommand{\fA}{{\mathfrak{A}}}

     \newcommand{\fP}{{\mathfrak{P}}}

    \newcommand{\ab}{{\mathrm{ab}}}

    \newcommand{\Aut}{{\mathrm{Aut}}}

    \newcommand{\Char}{{\mathrm{Char}}}

    \newcommand{\End}{{\mathrm{End}}}

    \newcommand{\Gal}{{\mathrm{Gal}}} \newcommand{\GL}{{\mathrm{GL}}}
    
    \newcommand{\Hom}{{\mathrm{Hom}}}
    
    \renewcommand{\Im}{{\mathrm{Im}}}
    \newcommand{\Ind}{{\mathrm{Ind}}}

    \newcommand{\ord}{{\mathrm{ord}}} \newcommand{\rk}{{\mathrm{rank}}}
     \newcommand{\Pic}{\mathrm{Pic}}

    \renewcommand{\mod}{\ \mathrm{mod}\ }

    \newcommand{\Sel}{{\mathrm{Sel}}}

    \newcommand{\SL}{{\mathrm{SL}}}

      \newcommand{\Supp}{{\mathrm{Supp}}}
    \newcommand{\tor}{{\mathrm{tor}}}

    \newcommand{\Vol}{{\mathrm{Vol}}}\newcommand{\vol}{{\mathrm{vol}}}

     \newcommand{\M}{\mathrm{M}}
        \newcommand{\Tr}{\mathrm{Tr}}

\newcommand{\matrixx}[4]{\begin{pmatrix}
#1 & #2 \\ #3 & #4
\end{pmatrix} }        

    \font\cyr=wncyr10

    \newcommand{\Sha}{\hbox{\cyr X}}
    \newcommand{\wh}{\widehat}

    \newcommand{\ov}{\overline}

    \newcommand{\lra}{\longrightarrow}
    \newcommand{\ra}{\rightarrow}

    \newcommand{\N}{\mathrm{N}}

                         \newcommand{\Red}{\mathrm{Red}}
                            \newcommand{\hilbert}[4]{\left(\frac{#1,#2}{#3;#4}\right)}

\newcommand{\lrb}[1]{\left(#1\right)}

\newcommand{\B}{\ensuremath{{\mathbb{B}}}}
\newcommand{\C}{\ensuremath{{\mathbb{C}}}}

\newcommand{\F}{\ensuremath{{\mathbb{F}}}}

 \newcommand{\BetaI}[1]{{\left\{#1 \right\} }}
\newcommand{\E}{\ensuremath{{\mathbb{E}}}}
\renewcommand{\L}{\ensuremath{{\mathbb{L}}}}
\newcommand{\Nm}{\ensuremath{\text{Nm}}}
\newcommand{\WaldsI}{I(\varphi,\chi)}
\newcommand{\zxz}[4]{\begin{pmatrix} #1 & #2 \\ #3 & #4 \end{pmatrix}}

\newcommand{\yh}[1]{}

\newcommand{\Cor}[1]{}
\newcommand{\New}[1]{#1}
\newcommand{\yhb}[1]{}

    \theoremstyle{plain}
    \newtheorem{thm}{Theorem}[section] \newtheorem{coro}[thm]{Corollary}
    \newtheorem{lem}[thm]{Lemma}  \newtheorem{prop}[thm]{Proposition}
    \newtheorem {conj}[thm]{Conjecture} \newtheorem{defn}[thm]{Definition}

\theoremstyle{remark} \newtheorem{remark}{Remark}[section]
\theoremstyle{remark} 
\theoremstyle{remark} 

    \numberwithin{equation}{section}

\begin{document}
\title[An explicit Gross-Zagier formula related to the Sylvester Conjecture]{An explicit Gross-Zagier formula related to the Sylvester Conjecture}
\author[Yueke Hu, Jie Shu and Hongbo Yin]{Yueke Hu, Jie Shu and Hongbo Yin}

\begin{abstract}
Let $p\equiv 4,7\mod 9$ be a rational prime number such that $3\mod p$ is not a cube. In this paper we prove the 3-part of $|\Sha(E_p)|\cdot |\Sha(E_{3p^2})|$ is as predicted by the Birch and Swinnerton-Dyer conjecture, where $E_p: x^3+y^3=p$ and $E_{3p^2}: x^3+y^3=3p^2$ are the elliptic curves related to the Sylvester conjecture and cube sum problems.
\end{abstract}
\address{Department of Mathematics, ETH, Zurich, Switzerland}
\email{huyueke2012@gmail.com}

\address{School of Mathematical Sciences, Tongji University, Shanghai 200092}
\email{shujie@tongji.edu.cn}
\address{School of Mathematics, Shandong University, Jinan 250100,
P.R.China}
\email{yhb2004@mail.sdu.edu.cn}

\thanks{Yueke Hu is supported by SNF-169247; Jie Shu is supported by NSFC-11701092; Hongbo Yin is partially supported by NSFC-11701548.}

\maketitle


\section{Introduction}
In this paper we are concerned about the explicit Gross-Zagier formula and the full BSD conjecture for the elliptic curves which are related to the Sylvester conjecture. The motivation comes from the cube sum problem. A nonzero rational number is called a cube sum if it is of the form $a^3+b^3$ with $a,b\in \BQ^\times$. For any $n\in \BQ^\times$, let $E_n$ be the elliptic curve over $\BQ$ defined by the projective equation $x^3+y^3=nz^3$ with the distinguished point $(1:-1:0)$.
If $n$ is not a cube or twice a cube of nonzero rationals, then $n$ is a cube sum if and only if $E_n(\BQ)$ has a point of infinite order.
A famous conjecture concerning the cube sums, attributed to Sylvester, is the following


\begin{conj}[Sylvester \cite{Sylv}, Selmer \cite{Selmer51}]\label{Sylv}
Any prime number $p\equiv 4,7,8\mod 9$ is a cube sum.
\end{conj}
For a good summary of this conjecture, please refer to \cite{DV1,DV17}. For an odd prime $p\geq 5$, a 3-descent \cite{Satge,DV1} gives that
\[\rk_\BZ E_p(\BQ)\leq\begin{cases}0,& p\equiv 2,5\mod 9;\\
                         1,&  p\equiv 4,7,8\mod 9;\\
                         2,&  p\equiv 1\mod 9.\end{cases}\]
Let $\epsilon(E_p)$ be the sign in the functional equation of the Hasse-Weil L-function $L(s,E_p)$. From \cite{Liverance}, we know 
\[\epsilon(E_p)=\begin{cases}
                         -1,&  p\equiv 4,7,8\mod 9;\\
                         +1,&  \text{otherwise}.\end{cases}\]
Then the Birch and Swinnerton-Dyer (BSD) conjecture implies the Sylvester conjecture.  In 1994,  Elkies announced a proof of Conjecture \ref{Sylv} for all primes $p\equiv 4,7\mod 9$, but without any detailed publication. However, Dasgupta and Voight \cite{DV17} proved the following weaker theorem using  method substantially different from that of Elkies.

\begin{thm}
Let $p\equiv 4,7\mod 9$ be a rational prime number such that $3\mod p$ is not a cubic residue. Then
$p$ and $p^2$ are cube sums. 
\end{thm}

Dasgupta and Voight proved the above theorem by establishing the non-triviality of certain related Heegner points. By the work of Gross-Zagier \cite{GZ1986} and Kolyvagin \cite{Kolyvagin1990}, the nontriviality of Heegner points implies that the rank part of the BSD conjecture for $E_p$ is true. 

If $\ell\nmid 6p$ is a prime, then $E_p$ has good reduction at $\ell$. Then Perrin-Riou \cite{PR1987} and Kobayashi \cite{Koba2013} proved that the $\ell$-part full BSD conjecture holds for $E_p$.  Since $E_p$ has potential good ordinary reduction at $p$,  the $p$-part full BSD conjecture of $E_p$ is also true by the work of Li-Liu-Tian \cite{LLT}. To summarize,  the following theorem is known. 
\begin{thm}Let $p\equiv 4,7\mod 9$ be a rational prime number such that $3\mod p$ is not a cubic residue. Then
\begin{enumerate}
\item[1.]  $\ord_{s=1}L(s,E_p)=\rk_{\BZ} E_p(\BQ)=1$;
\item[2.] The Tate-Shafarevich group $\Sha(E_p)$ is finite, and for any prime $\ell\nmid 6$, the $\ell$-part of $|\Sha(E_p)|$ is as predicted by the Birch and Swinnerton-Dyer conjecture for $E_p$.
\end{enumerate}
\end{thm}

But for the primes $\ell=2,3$, there are no results known for the $\ell$-part full BSD conjecture of $E_p$. In this paper, we adopt a similar method as in \cite{CST17} to approach the 3-part full BSD conjecture of $E_p$ and $E_{3p^2}$ by comparing to  an explicit Gross-Zagier formula.

Let $\Sha(E_p)$ denote the Shafarevich-Tate group of $E_p$, $E_p(\BQ)_\tor$ denote the torsion group of $E_p(\BQ)$, $\Omega_p$ denote the minimal real period of $E_p$, $\wh{h}(\cdot)$ denote the N\'eron-Tate height of $E_p$ over $\BQ$ and $c_\ell$ denote the Tamagawa number of $E_p$ at a prime $\ell$. From \cite{DV17}, we know that $E_p(\BQ)$ resp. $E_{3p^2}(\BQ)$ has rank $1$ resp. $0$. Let $P$ be a generator of the free part of $E_p(\BQ)$.  Then the BSD conjecture predicts that
\begin{equation*}
 |\Sha(E_p)|=\frac{L'(1,E_p)}{\Omega_p\cdot \widehat{h}_\BQ(P)}\cdot \frac{ |{E_p}(\BQ)_\tor|^2}{\prod_\ell c_\ell(E_p)},
\end{equation*}
where $\ell$ runs through all prime numbers. Similarly,for $E_{3p^2}(\BQ)$ the BSD conjecture predicts that
\begin{equation*}
\left|\Sha(E_{3p^2})\right|=\frac{L(1,E_{3p^2})}{\Omega_{3p^2} }\cdot \frac{ |{E_{3p^2}}(\BQ)_\tor|^2}{\prod_\ell c_\ell(E_{3p^2})}.
\end{equation*}
Combining these two formulae, we shall expect that
\begin{equation}\label{bsd1}
|\Sha(E_p)|\cdot|\Sha(E_{3p^2})|=\frac{L'(1,E_p)}{\Omega_p\cdot \widehat{h}_\BQ(P)}\cdot\frac{L(1,E_{3p^2})}{\Omega_{3p^2}}\cdot \frac{ |{E_p}(\BQ)_\tor|^2}{\prod_\ell c_\ell(E_p)}\cdot \frac{ |{E_{3p^2}}(\BQ)_\tor|^2}{\prod_\ell c_\ell(E_{3p^2})}.
\end{equation}
Our main result is the following.
\begin{thm}\label{Main}
Let $p\equiv 4,7\mod 9$ be a rational prime number such that $3\mod p$ is not a cubic residue. Then the both sides of (\ref{bsd1}) are nonzero rational numbers and the exponents of $3$ in both sides of (\ref{bsd1}) are equal as expected.\end{thm}

In the following, we sketch the proof of Theorem \ref{Main}. Let $\CH=\{z\in \BC:\Im(z)>0\}$ be the Poinc\'are upper half plane and $SL_2(\BZ)$ acts on $\CH$ by fractional linear transformations.   Let $\Gamma_0(3^5)\subset \SL_2(\BZ)$ be the congruence group of level $3^5$ which consists of matrices 
\[\matrixx{a}{b}{c}{d}\textrm{ with }c\equiv 0\mod 3^5.\] 
Then $Y_0(3^5)=\Gamma_0(3^5)\backslash \CH$  is an affine smooth curve over $\BQ$ and let $X_0(3^5)$ be its projective closure.  

Fix $K=\BQ(\sqrt{-3})\subset \BC$ with $\CO_K=\BZ[\omega]$ its ring of integers, where $\omega=\frac{-1+\sqrt{-3}}{2}$. We carefully embed $K$ into $\M_2(\BQ)$ as follows. Once such an embedding is given, the group $K^\times $ of invertible elements acts on $\CH$ through fractional linear transformations and there is a unique point in $\CH$ which is invariant under the action of $K^\times$. There are exactly two embeddings $\rho:K\hookrightarrow \M_2(\BQ)$ with fixed point $\tau=(2p\omega-9)/(9p\omega-36)\in \CH$ and we choose the normalized one, i.e. we have
\[\rho(t)\begin{pmatrix}\tau\\1\end{pmatrix}=t\begin{pmatrix}\tau\\1\end{pmatrix},\quad \textrm{for any $t\in K$}.\]

For any $n\in \BQ^\times$, the elliptic curve $E_n$ has Weierstrass equation $$y^2=x^3-432n^2$$ and has complex multiplication by $\CO_K$ over $K$. We fix the complex multiplication $[\ ]:\CO_K\simeq \End_{\ov{\BQ}}(E_n)$ by $[\omega](x,y)=(\omega x,y)$.  The image of  $\tau\in \CH$ defines a CM point on $X_0(3^5)$. Let $f:X_0(3^5)\ra E_9$ be the the natural modular parametrization. The Heegner point $f(\tau)$ is defined over the ring class field $H_{9p}$ over $K$ of conductor $9p$. 

Let $p\equiv 4,7\mod 9$ be a prime.  Let 
$\chi:\Gal(\ov{K}/K)\ra \CO_K^\times$ be the character given by $\chi(\sigma)=(\sqrt[3]{3p})^{\sigma-1}$. 
 Define the Heegner cycle 
$$P_\chi(f)=\sum_{\sigma\in \Gal(H_{9p}/K)}f(\tau)^\sigma\otimes \chi(\sigma)\in E_9(H_{9p})\otimes_\BQ K.$$ 
The base change L-function $L(s,E_9,\chi)$ has  sign $-1$ and has a decomposition
\[L(s,E_9,\chi)=L(s,E_p)\cdot L(s,E_{3p^2}).\]
By the work \cite{DV17}, we know that $L(s,E_p)$ has a zero of order $1$ at $s=1$ and $L(1,E_{3p^2})\neq 0$.
The morphism $f$  is a test vector for  the pair $(E_9,\chi)$, i.e. there is a nontrivial relation between the central value of the derivative of the L-function $L(s,E_9,\chi)$ and the height of  the Heegner cycle $P_\chi(f)$. More precisely, we have the following result(see Theorem \ref{thm:GZ}). 
\begin{thm} \label{thm:GZ1}
For primes $p\equiv 4,7 \mod 9$,  we  have the following explicit height formula of Heegner cycles:
\[\frac{L'(1,E_{p})L(1,E_{3p^2})}{\Omega_{p}\Omega_{3p^2}}=2^{\alpha}\cdot \left\langle P_{\chi}(f),P_{\chi^{-1}}(f)\right\rangle_{K,K}\]
where $\alpha=0$ if $p\equiv 4\mod 9$ and $\alpha=-1$ if $p\equiv 7\mod 9$, and $\langle \cdot,\cdot\rangle_{K,K}$ denotes the $K$-linear N\'eron-Tate height pairing of $E_9$ over $K$
\end{thm}

For the definition of $\langle \cdot,\cdot\rangle_{K,K}$, see for example \cite[page 2531]{CST14}. We remark that we don't need the hypothesis that $3$ is not a cube modulo $p$ in the above theorem.

Comparing this explicit Gross-Zagier fomula with the product formula (\ref{bsd1})  of  full BSD conjectures for $E_p$ and $E_{3p^2}$, Theorem \ref{Main} follows from the $\sqrt{-3}$-nondivisibility of Heegner points.

This paper is organized as follows. In Section \ref{Heegner}, we give the construction of the Heegner points and study the Galois actions on the Heegner points via modular actions. In Section \ref{Non-divisibility}, we briefly recall the non-triviality of the Heegner points from \cite{DV17} and study the 3-nondivisibility of the Heegner points. In Section \ref{E-G-Z}, we establish the explicit Gross-Zagier formula for the Heegner points (Theorem \ref{thm:GZ1}) assuming a local period integral result. In Section \ref{E-B-SD}, we prove Theorem \ref{Main} by comparing the explicit Gross-Zagier formula and the full BSD conjecture. 
In Section \ref{minimal vector} we briefly review the compact induction theory from \cite{BushnellHenniart:06a} and the results on local Waldspurger's period integral using minimal vectors from \cite{HN18}. In Section \ref{newforms} we establish explicit relation between minimal vectors and newforms, and derive special cases of the local Waldspurger's period integral for newforms used in section \ref{E-B-SD}. We also give more general results for newforms which are of independent interest and can be useful for further applications.

\noindent {\bf Acknowledgements.} The authors would like to thank Professor Ye Tian for useful conversations and  encouragement. We also would like to thank John Voight and Samit Dasgupta who send their unpublished work and provide many useful discussions and helps. We also thank Li Cai, who provide many discussions about the matrix coefficients of supercuspidal representations and Jianing Li who provides many helps on the SageMath systems. We would like to thank the referee for helpful advices which motivate us to improve the treatments for 3-nondivisibility of the Heegner points and local period integrals in the current version.


\section{Modular Actions on Heegner Points}\label{Heegner}
\subsection{The Modular Curves and Modular Actions}
Let  $X$ be an algebraic curve defined over $\BQ$ and  $F$  a field extension of $\BQ$.  Denote by $\Aut_F(X)$ the group of algebraic automorphisms of $X$ which are defined over $F$.
Let $$\CH=\{z\in \BC|\, \Im(z)>0\}$$ be the Poinc\'are upper half plane. The group $\GL_2(\BQ)^+$ acts on $\CH$ by linear fractional transformations.

Let $U_0(3^5)$ be the open compact subgroup of $\GL_2(\widehat{\BZ})$ consisting of matrices
$\left (
\begin{array}{cc}
a&b\\c&d
\end{array}
\right )$ such that $c\equiv 0\mod 3^5$, and let
$\Gamma_0(3^5)=\GL_2(\BQ)^+\cap U_0(3^5)$. Let $X_0(3^5)$ be the
modular curve over $\BQ$ of level $\Gamma_0(3^5)$ whose underlying Riemann surface is
\[X_0(3^5)(\BC)=\GL_2(\BQ)^+\backslash \left(\CH\sqcup\BP^1(\BQ)\right )\times \GL_2(\BA_f)/U_0(3^5)\simeq \left( \Gamma_0(3^5)\backslash \CH\right)\sqcup \left(\Gamma_0(3^5)\backslash \BP^1(\BQ)\right),\]
where $\BA_f$ denote the finite ad\`ele of $\BQ$.
Define $N$ to be the normalizer of $\Gamma_0(3^5)$ in $\GL_2^+(\BQ)$.
It follows from \cite[Theorem 1]{KM1988} that the linear fractional transformation action of $N$ on $X_0(3^5)$ induces an isomorphism
\[N/\BQ^\times\Gamma_0(3^5)\simeq \Aut_{\ov{\BQ}}(X_0(3^5)).\]
Moreover, all the algebraic automorphisms in $\Aut_{\ov{\BQ}}(X_0(3^5))$ are defined over $K$.
We identify $\Aut_{\ov{\BQ}}(X_0(3^5))$ with $N/\BQ^\times\Gamma_0(3^5)$ by this isomorphism. By \cite[Theorem 8]{AL1970}, and \cite{Ogg80}, the quotient group $N/\BQ^\times\Gamma_0(3^5)\simeq S_3\rtimes \BZ/3\BZ$, where $S_3$ denotes the symmetric group with $3$ letters which is generated by the Atkin-Lehner operator $W=\begin{pmatrix}0&1\\-3^5&0\end{pmatrix}$ and the matrix $A=\begin{pmatrix}28&1/3\\3^4&1\end{pmatrix}$, and the subgroup $\BZ/3\BZ$ is generated by the matrix $B=\begin{pmatrix}1&0\\3^4&1\end{pmatrix}$.

Put
\[U=\langle U_0(3^5),W,A\rangle\subset \GL_2(\BA_f).\]
Then $\BQ^\times\backslash \BQ^\times U$ is an open compact subgroup of $\BQ^\times\backslash \GL_2(\BA_f)$.
Put $$\Gamma =\GL_2(\BQ)^+\cap U=\langle \Gamma_0(3^5),W,A\rangle,$$ and let $X_\Gamma$ be the
modular curve over $\BQ$ of level $\Gamma$ whose underlying Riemann surface is
\[X_\Gamma(\BC)=\GL_2(\BQ)^+\backslash\left (\CH\sqcup\BP^1(\BQ)\right )\times \GL_2(\BA_f)/U\simeq (\Gamma \backslash \CH)\sqcup (\Gamma\backslash \BP^1(\BQ)).\]
Then $X_\Gamma$ is a smooth projective curve over $\BQ$ of genus $1$, and $X_\Gamma$ has three cusps
\[\Gamma\backslash \BP^1(\BQ)=\{[\infty],[1/9],[2/9]\}.\] The cusp $[\infty]$ is rational over $\BQ$, and the cusps $[1/9]$ and $[2/9]$ are both defined over $K$. We identify $X_\Gamma$ with an elliptic curve over $\BQ$ with $[\infty]$ as its zero element. Let $N_\Gamma$ be the normalizer of $\Gamma$ in $\GL_2(\BQ)^+$. Then we have a natural embedding
\[\Phi: N_\Gamma/ \BQ^\times\Gamma\hookrightarrow \Aut_{\ov{\BQ}}(X_\Gamma)\simeq \CO_K^\times\ltimes X_\Gamma(\ov{\BQ}),\]
where $\CO_K^\times $ embeds into $\Aut_{\ov{\BQ}}(X_\Gamma)$ by complex multiplications and $X_\Gamma(\ov{\BQ})$ embeds into $\Aut_{\ov{\BQ}}(X_\Gamma)$ by translations.
The matrices
\[B=\begin{pmatrix}1&0\\3^4&1\end{pmatrix},\quad C=\begin{pmatrix}1&1/9\\-3^3&-2\end{pmatrix}\]
lie in $N_\Gamma$, and hence induce automorphisms of $X_\Gamma$.

The elliptic curves $E_n$ are all endowed with complex multiplication by $K$ and we fix the complex multiplication $[\cdot]:\CO_K\simeq \End_K(E_n)$ by $[-\omega](x,y)=(\omega x,-y)$. We will always take the simple Weierstrass equation $y^2=x^3-2^4\cdot 3$ for the elliptic curve $E_9$, unless stated otherwise. 
\begin{prop}\label{modular-curve}
\begin{itemize}
\item[1.] The elliptic curve $(X_\Gamma,[\infty])$ is isomorphic to $E_9$ over $\BQ$.
\item[2.] We have an embedding
\[\Phi:N_\Gamma/ \BQ^\times\Gamma\hookrightarrow \CO_K^\times \ltimes (\Gamma\backslash \BP^1(\BQ))\subset \Aut_{\ov{\BQ}}(X_\Gamma).\]
Moreover,  for any point $P\in X_\Gamma$, we have
\[\Phi(B)(P)=[\omega^2]P,\quad \Phi(C)(P)=[\omega^2]P+[1/9].\]
 In particular, the automorphisms $\Phi(B)$ and $\Phi(C)$ are defined over $K$.
\end{itemize}
\end{prop}
Note that there exists a unique isomorphism $X_\Gamma\ra E_9$ over $\BQ$ such that the cusp $[1/9]$ has coordinates $(0,4\sqrt{-3})$. We use this isomorphism to identify $X_\Gamma$ with $E_9$.
\begin{proof}
It is known from \cite{DV17} that $E_9$ is the natural quotient of $X_0(3^5)$ by the finite group $S_3$.
Since the automorphism group of the elliptic curve $E_9$ is isomorphic to $ \CO_K^\times$, we have \[\Aut_{\ov{\BQ}}(X_\Gamma)\simeq \CO_K^\times\ltimes X_\Gamma(\ov{\BQ}).\]
Then for any $M\in N_\Gamma$ and $P\in X_\Gamma$, $\Phi(M)(P)=[\alpha]P+S$ for some $\alpha\in \CO_K^\times,S\in X_\Gamma(\ov{\BQ})$. Taking $P=[\infty]$, we see $S=\Phi(M)([\infty])\in \Gamma\backslash \BP^1(\BQ)$. The formulae for $\Phi(B)$ and $\Phi(C)$ are taken from \cite{DV17}, which can also be verified numerically using SageMath.
\end{proof}

Let $V\subset U_0(3^5)$  be the subgroup consisting of matrices
$\left (\begin{array}{cc}a&b\\c&d
\end{array} \right )$ with $a\equiv d\mod 3$, and put $U_0=\langle V,W,A \rangle$.  Let $X_\Gamma^0$ be the modular curve over $\BQ$ whose underlying Riemann
surface is
\[X_\Gamma^0(\BC)=\GL_2(\BQ)^+\backslash\left (\CH\bigsqcup\BP^1(\BQ)\right )\times \GL_2(\BA_f)/U_0.\]
By class field theory, $ \BQ^\times_+\wh{\BZ}^\times/\BQ^\times_+\det(U_0)\simeq \Gal(K/\BQ).$
Noting that $\GL_2(\BQ)^+\cap U_0=\Gamma$, we see that the
modular curve $X_\Gamma^0$ is isomorphic to $X_\Gamma\times_\BQ K$ as a curve over $\BQ$ (cf. \cite[Chapter 6]{Shimurabook}). Usually, we denote by $[z,g]_{U_0}$  the point  on  $X_\Gamma^0$ which is represented by the pair $(z,g)$ where $z\in \CH$ and $g\in \GL_2(\BA_f)$. The curve $X_\Gamma^0$ is not geometrically connected
and has two connected components over $\BC$. Put
$$U/U_0=\langle\epsilon\rangle,\quad \epsilon=\left (\begin{array}{cc}1&0\\0&-1\end{array}\right ).$$
The non-trivial Galois action of $\Gal(K/\BQ)$ on $X_\Gamma^0$ is given by
the right translation of $\epsilon$ on $X_\Gamma^0$. We have
\[\Aut_\BQ(X_\Gamma^0)=\Aut_K(X_\Gamma)\rtimes\Gal(K/\BQ)\simeq (X_\Gamma(K)\rtimes \CO_K^\times)\rtimes \Gal(K/\BQ).\]
Let $N_{\GL_2(\BA_f)}(U_0)$ be the normalizer of $U_0$ in
$\GL_2(\BA_f)$. Then there is a natural homomorphism
$$ N_{\GL_2(\BA_f)}(U_0)/U_0\longrightarrow \Aut_\BQ(X_\Gamma^0)$$ induced by right
translation on $X_\Gamma^0$: for $P=[z,g]_{U_0}\in X_\Gamma^0$ and $x\in N_{\GL_2(\BA_f)}(U_0)$
\[P\mapsto P^x=[z,gx]_{U_0}.\] 
An element $g\in
N_{\GL_2(\BA_f)}(U_0)$ maps one component of $X_\Gamma^0$ onto the other if and
only if it has image $-1$ under the composition of the following
morphisms:
\[\xymatrix{\GL_2(\BA_f)=\GL_2(\BQ)^+\GL_2(\wh{\BZ})\ar[r]^{\quad\quad \quad \quad \quad\det}&\BQ^{\times}_+\widehat{\BZ}^\times\ar[r]&\BZ_3^\times/(1+3\BZ_3),}\]
where $\wh{\BZ}=\prod_\ell \BZ_\ell$ and the last morphism is trivial on $\BQ^\times_+$ and, on $\wh{\BZ}^\times$, it is the projection from
$\widehat{\BZ}^\times$ to its $3$-adic factor composed with $\mod 3$.

\subsection{The Modular Actions on Heegner Points}\label{heegner}
Let $p\equiv 4,7\mod 9$ be a rational prime number. Let $\rho:K\rightarrow \M_2(\BQ)$ be the normalised embedding with fixed point $\tau=(2p\omega-9)/(9p\omega-36)\in \CH$, i.e. we have
\[\rho(t)\begin{pmatrix}\tau\\1\end{pmatrix}=t\begin{pmatrix}\tau\\1\end{pmatrix},\quad \textrm{for any $t\in K$}.\]
Here it's matrix multiplication on the left hand side and scalar multiplication on the right hand side.
Note that $$\tau=M\omega,\quad M=\matrixx{2}{-1}{9}{-4}\matrixx{\frac{p}{9}}{0}{0}{1}.$$
Then the embedding $\rho:K\rightarrow\M_2(\BQ)$ is explicitly given by
\[\rho(\omega)=M\matrixx{-1}{-1}{1}{0}M^{-1}=
\begin{pmatrix}2p+8+36/p&-4p/9-2-9/p\\9p+36+144/p&-2p-9-36/p\end{pmatrix}.\]
Let $R_0(3^5)$ be the standard Eichler order of discriminant $3^5$ in $\M_2(\BQ)$.  For any integer $c\geq 1$, let $\CO_c$ be the order of $K$ of conductor $c$ and let $H_c$ be the ring class field of conductor $c$.  
Then $K\cap R_0(3^5)=\CO_{9p}$. Let $\CO_{K,3}$ be the completion of $\CO_K$ at the unique place above $3$. We have
\[ \CO_{K,3}^\times/\BZ_3^\times(1+9\CO_{K,3})=\langle \omega_3\rangle\times\langle1+3\omega_3\rangle\cong \BZ/3\BZ\times\BZ/3\BZ,\]
where $\omega_3$ is the image of $\omega$ into $\CO_{K,3}^\times$. Considered as elements in $GL_2(\BA_f)$ with other components 1, it is straightforward to verify that $\omega_3$ and $1+3\omega_3$ normalize $U_0$, and hence we have an embedding
\[\CO_{K,3}^\times/\BZ_3^\times(1+9\CO_{K,3})\hookrightarrow \Aut_\BQ(X_\Gamma^0).\]
If $p\equiv 7\mod 9$, it is straight-forward to verify that  the element
\[
w=M\matrixx{1}{1}{0}{-1}M^{-1}=\begin{pmatrix}-2p-17&4p/9+4\\-9p-72&2p+17\end{pmatrix}.
\]
is a nontrivial normalizer of $K^\times$ in $\GL_2(\BQ)$ and $w\epsilon$ normalizes $U_0$, and hence $w\epsilon$ also induces an automorphism of $X_\Gamma^0$.

\begin{thm}\label{modular-action}
\begin{itemize}
\item[1.] For any point $P\in X_\Gamma^0$, we have
\[P^{1+3\omega_3}=[\omega^2]P,\]
and
\[P^{\omega_3}=\left\{\begin{aligned}{ [\omega^2]P+(0,4\sqrt{-3})},&\quad p\equiv 4\mod 9 ;\\{[\omega]P+(0,4\sqrt{-3})},&\quad p\equiv 7\mod 9 . \end{aligned}\right.\]
\item[2.]  Suppose $p\equiv 7\mod 9$. For any point $P\in X_\Gamma^0$, we have
\[P^{w\epsilon}=[\omega^{\frac{p-7}{9}}]P-(0,4\sqrt{-3}).\]
\end{itemize}
\end{thm}
\begin{proof}
Since $\omega_3$, $1+3\omega_3$  and $w\epsilon$ all have determinant $\equiv 1\mod 3$, when identified as elements in $\Aut_\BQ(X_\Gamma^0)$ they actually lie in the subgroup $\Aut_K(X_\Gamma)$.
Let $P=[z,1]_{U_0}$ for $z\in \CH$ be a point on $X_\Gamma^0$. We have
\[B(1+3\omega_3)A^2=\left(\begin{pmatrix}60p+837/p+214&2p/3+9/p+7/3\\ 5130p+71145/p+18252&57p+765/p+199\end{pmatrix}_3, BA^2\right)\in V,\]
where the subscript $3$ denotes the $3$-adic component of the adelic matrices.
Then
\[P^{1+3\omega_3}=[z,1+3\omega_3]_{U_0}=[B(z),B(1+3\omega_3)]_{U_0}=[B(z),1]_{U_0}=\Phi(B)P=[\omega^2]P.\]

If $p\equiv 4\mod 9$, then
\[C\omega_3 A^2=\left(\begin{pmatrix}     867p + 11635/p + 2685    &   31p/3 + 416/3p + 32\\
-20808p - 281925/p - 64440   &    -248p - 3360/p - 768\end{pmatrix}_3, CA^2\right)\in V.\]
and hence
\[P^{\omega_3}=\Phi(C)(P)=[\omega^2]P+(0,4\sqrt{-3}).\]
If $p\equiv 7\mod 9$, then
\[BC\omega_3 A^2=\left(\begin{pmatrix}      867p + 11635/p + 2685    &   31p/3 + 416/3p + 32\\
49419p + 660510/p + 153045   &    589p + 7872/p + 1824\end{pmatrix}_3, BCA^2\right)\in V.\]
and hence
\[P^{\omega_3}=\Phi(BC)(P)=[\omega]P+(0,4\sqrt{-3}).\]

Suppose $p\equiv 7\mod 9$. It can verified that
\[\left\{\begin{aligned}BC^2w\epsilon A^2 \in  V,&\quad \frac{p-7}{9}\equiv 0\mod 3; \\
C^2w\epsilon A^2 \in  V ,&\quad \frac{p-7}{9}\equiv 1\mod 3;\\ B^2C^2w\epsilon A^2 \in  V ,&\quad \frac{p-7}{9}\equiv 2\mod 3. \end{aligned}\right.\]
Hence the second assertion follows.

\end{proof}

Let $\sigma: \wh{K}^\times\ra \Gal(K^\ab/K)$ be the Artin reciprocity law and we denote by $\sigma_t$ the image of $t\in\wh{K}^\times$.   Let $P_0=[\tau,1]$ be the CM point on $X_\Gamma^0$.
\begin{thm}\label{thm:Galois}
\begin{itemize}
\item[1.] The point $P_0\in X_\Gamma^0(H_{9p})$ satisfies
\[P_0^{\sigma_{1+3\omega_3}}=[\omega^2]P_0,\]
and
\[P_0^{\sigma_{\omega_3}}=\left\{\begin{aligned} {[\omega^2]P_0+(0,4\sqrt{-3})},&\quad p\equiv 4\mod 9 \\{[\omega]P_0+(0,4\sqrt{-3})},&\quad p\equiv 7\mod 9  \end{aligned}\right.\]
\item[2.] Suppose $p\equiv 7\mod 9$. We have
\[\ov{P_0}=[\omega^{\frac{p-7}{9}}]P_0-(0,4\sqrt{-3}).\]
\end{itemize}
\end{thm}
\begin{proof}
By Shimura's reciprocity law \cite[Theorems 6.31 and  6.38]{Shimurabook}, we have
\[P_0^{\sigma_t}=P_0^t=[\tau, t],\quad t\in \wh{K}^\times.\]
Since $\wh{K}^\times\cap U_0=\wh{\CO_{9p}}^\times,$ we see $P_0$ is defined over the ring class field $H_{9p}$, and the Galois actions of $\sigma_{\omega_3}$ and $\sigma_{1+3\omega_3}$ are clear from Theorem \ref{modular-action}.
\end{proof}
\begin{prop}\label{LCF}
\begin{itemize}
\item[1.] We have $H_{9p}=H_{3p}(\sqrt[3]{3})$ with $\Gal(H_{9p}/H_{3p})= \langle\sigma_{1+3\omega_3}\rangle\simeq {\BZ/3\BZ}$, and
\[\left(\sqrt[3]{3}\right)^{\sigma_{1+3\omega_3}-1}=\omega^2.\]
\item[2.] We have $\left(\sqrt[3]{3}\right)^{\sigma_{\omega_3}-1}=1$ and
\[\left(\sqrt[3]{p}\right)^{\sigma_{\omega_3}-1}=\left\{\begin{aligned} {\omega^2},&\quad p\equiv 4\mod 9 \\{\omega},&\quad p\equiv 7\mod 9  \end{aligned}\right.\]
\end{itemize}
\end{prop}
\begin{proof}
For any place $w$ of $K$, let $K_w$ denote the completion of $K$ at the place $w$, and let $\hilbert{\cdot}{\cdot}{K_w}{3} $ be the third Hilbert symbol over $K_w$, see for example \cite[Chapter V, Section 3]{Neukirchbook1}. We have the decomposition of ideal $7\CO_K=(1+3\omega)(1+3\omega^2)$. Let $v$ be the place corresponding to the prime ideal $(1+3\omega)$. Then
\[\left(\sqrt[3]{3}\right)^{\sigma_{1+3\omega_3}-1}=\hilbert{1+3\omega_3}{ 3}{K_3}{3}.\]
It is an important property of Hilbert symbol that
\[\prod_w \hilbert{1+3\omega_w}{ 3}{K_w}{3}=1,\]
where $\omega_w$ denotes the image of $\omega$ in $K_w$ and $w$ runs through all places of $K$. Since the symbol is trivial whenever $w\neq 3, v$, we have by \cite[Chapter V, Proposition 3.4]{Neukirchbook1} that 
\[\left(\sqrt[3]{3}\right)^{\sigma_{1+3\omega_3}-1}=\hilbert{1+3\omega_3}{ 3}{K_3}{3}=\hilbert{1+3\omega_v}{ 3}{K_v}{3}^{-1}=3^{-2}\mod (1+3\omega)=\omega^2.\]

Since $p\equiv 1\mod 3$, the prime $p$ splits in $K$ and let $v$ and $\ov{v}$ be the two places of $K$ above $p$. Then similarly
\[\left(\sqrt[3]{p}\right)^{\sigma_{\omega_3}-1}=\hilbert{\omega_3}{ p}{K_3}{3}=\hilbert{\omega_v}{ p}{K_v}{3}^{-1}\cdot \hilbert{\omega_{\ov{v}}}{ p}{K_{\ov{v}}} {3}^{-1}=\omega^{-\frac{p-1}{3}}.\]
\end{proof}

The elliptic curve $E_1$ has Weierstrass equation $y^2=x^3-432$. Consider the isomorphism
\[\phi:E_9\lra E_1,\quad (x,y)\mapsto ((\sqrt[3]{3})^2x, 3y).\]
We have the following commutative diagram:
$$\xymatrix{E_9(H_{9p})^{\sigma_{1+3\omega_3}=\omega^2}\ar[d]^\phi\ar[rr]^{\Tr_{H_{9p}/L_{(3,p)}}}&&E_9(L_{(3,p)})^{\sigma_{1+3\omega_3}=\omega^2}\ar[d]^{\phi}\\
            E_1(H_{3p})\ar[rr]^{\Tr_{H_{3p}/L_{(p)}}}&&E_1(L_{(p)})}$$
where the field extension diagram is as follows:
\[\xymatrix@R=11pt{&H_{9p}=H_{3p}(\sqrt[3]{3})\ar@{-}[dl]^{3}\ar@{-}[dr]\ar@{-}[dd]&\\
            H_{3p}\ar@{-}[dd]_{(p-1)/3}&&H_9\ar@{=}[dd]\\
            &L_{(3,p)}=K(\sqrt[3]{3},\sqrt[3]{p})\ar@{-}[dr]^{3}\ar@{-}[d]^{3}\ar@{-}[dl]_{3}&\\
            L_{(p)}=K(\sqrt[3]{p})\ar@{-}[dr]^{3}&L_{(3p)}=K(\sqrt[3]{3p})\ar@{-}[d]^{3}&L_{(3)}=K(\sqrt[3]{3})\ar@{-}[dl]_{3}\\
            &K\ar@{-}[d]^{2}&\\
            &\BQ.&\\
            }\]
Put $Q=\phi(P_0)$ and $R=\Tr_{H_{3p}/L_{(p)}} Q$.
\begin{coro}\label{co: Galois}
The point $R\in E_1(L_{(p)})$ and satisfies
\[R^{\sigma_{\omega_3}}=\left\{\begin{aligned} {[\omega^2]R+(0,12\sqrt{-3})},&\quad p\equiv 4\mod 9, \\{[\omega]R+(0,-12\sqrt{-3})},&\quad p\equiv 7\mod 9,  \end{aligned}\right.\]
and if $p\equiv 7\mod 9$,
\[\ov{R}=[\omega^{\frac{p-7}{9}}]R+(0,12\sqrt{-3}).\]
\end{coro}
\begin{proof}
This is a consequence of Theorem \ref{thm:Galois}, Proposition \ref{LCF} and that 
\[\Tr_{H_{3p}/L_{(p)}} (0,12\sqrt{-3})=\frac{p-1}{3}(0,12\sqrt{-3}).\]
Recall that the cusp $(0,12\sqrt{-3})$ is a 3-torsion.
\end{proof}

\section{Non-divisibility of Heegner points}\label{Non-divisibility}
By the assumption that $3\mod p$ is not a cubic residue, we decompose $p\CO_{K}=\fp\ov{\fp}$ so that
\begin{equation}\label{h}
3^{\frac{p-1}{3}}\equiv \omega\mod \fp, \quad 3^{\frac{p-1}{3}}\equiv \omega^2\mod \ov{\fp}.
\end{equation}
Since $H_{3p}/K$ is totally ramified at $\fp$ and $\ov{\fp}$, let $\fP$ and $\ov{\fP}$ be the primes of $H_{3p}$ above $\fp$ and $\ov{\fp}$ respectively. We have
\[\CO_K/p\CO_K\cong \CO_K/\fp\bigoplus\CO_K/\ov{\fp}\subset \CO_{H_{3p}}/\fP\bigoplus \CO_{H_{3p}}/\ov{\fP}.\]

The main result (Proposition 5.2.8) in \cite{DV17} state that if $3$ is not a cube modulo $p$, the reduction
\[(R\mod \fP,R\mod \ov{\fP})\in E_1(\BF_p)^2\]
is not equal to the image of any torsion point in $E_1(L_{(p)})$, and hence the Heegner point $R$ is of infinite order. 

We sketch their strategy briefly.  Note the CM points $P_0$ considered in this paper are exactly those considered in \cite{DV17}. We will use the explicit coordinates of $R$ modulo $p$ later so we record them in the following lemma.
\begin{lem}\label{lem:coord-R}
Under the Weierstrass equation $y^2=x^3-432$ for $E_1$, 
\[
R\equiv\left\{ \begin{aligned} (12\omega^{i+2},-36) &\mod \fP,\\ (12\omega^{i+1},-36)&\mod \ov{\fP},\end{aligned}\right.
\]
where $i=0,1,2$ depends on $p$.
\end{lem}
\begin{proof}

As in the proof of \cite[Proposition 5.2.8]{DV17}, 
\begin{equation}\label{heegner-cong}
R=\Tr_{H_{3p/L_{(p)}}} Q\equiv \frac{p-1}{3}Q\equiv \begin{cases}Q,& p\equiv 4\mod 9,\\ -Q, & p\equiv 7\mod 9.\end{cases}
\end{equation}
By \cite[Proposition 5.2.1 and Lemma 5.2.4]{DV17}, the x-coordinate $x(Q)$ is $p$-adic integral and satisfies the following congruence (see \cite[(5.2.6)]{DV17})
\begin{equation}\label{xcoor1}
x(Q)\equiv x(Q)^p\equiv 12^p\omega^i(-3)^{(p-1)/6}\mod p\ov{\BZ}
\end{equation}
where $\ov{\BZ}$ denotes the ring of integral algebraic numbers and $i=0,1,2$ depends on $p$ (see \cite[Lemma 5.1.3]{DV17}).
Note that Dasgupta and Voight \cite{DV17} used the Weierstrass equation $y^2+y=3x^3-1$ for $E_1$, and we have adapted the coordinates to the equation $y^2=x^3-432$.
Since $-3\in \BF_p^\times$ is a square, $(-3)^{\frac{p-1}{6}}$ is a third root of unity modulo $p$. Note also that $\frac{p-1}{3}$ is always even in our case. Therefore, by (\ref{h}) and (\ref{xcoor1}), we conclude that
\begin{equation}\label{coordinates}
x(Q)\equiv\left\{ \begin{aligned} 12\omega^{i+2}&\mod \fP,\\12\omega^{i+1}&\mod \ov{\fP}.\end{aligned}\right.
\end{equation}
By Theorem \ref{thm:Galois}, we have
\[[1-\omega^\alpha]Q\equiv (0,12\sqrt{-3}) \mod \fP\quad  \textrm{(resp. $\ov{\fP}$)},\]
where $\alpha=2$ if $p\equiv 4\mod 9$ and $\alpha=1$ if $p\equiv 7\mod 9$. This implies that $(Q \mod \fP)$ and $(Q \mod \ov{\fP})$ both belongs to $ E_1(\BF_p)[3]$. It follows  that
\begin{equation}\label{coor3}
Q\equiv\left\{ \begin{aligned} (12\omega^{i+2},(-1)^{\alpha-1}\cdot 36) &\mod \fP,\\ (12\omega^{i+1},(-1)^{\alpha-1}\cdot 36)&\mod \ov{\fP}.\end{aligned}\right.
\end{equation}
Then the Lemma follows from (\ref{heegner-cong}) and (\ref{coor3}).
\end{proof}

Consider the reduction map
\[\mathrm{Red}:E_1(L_{(p)})\backslash \{O\} \lra \CO_{H_{3p}}/\fP\bigoplus \CO_{H_{3p}}/\ov{\fP},\quad T\mapsto (x(T)\mod \fP,x(T)\mod \ov{\fP}).\]
By \cite[Lemma 5.2.9]{DV17}, we have $E_1(L_{(p)})_\tor=E_1(K)_\tor$.
Let $D$ be the image of $E_1(L_{(p)})_\tor\backslash \{O\}$ under the reduction map. Then
\[D=\{(0,0),(12\omega^i,12\omega^i)_{i=0,1,2}\}\subset \CO_K/\fp\bigoplus\CO_K/\ov{\fp}.\]
On the other hand, by Lemma \ref{lem:coord-R}, the reduction $\Red(R)$ is not trivial and also doesn't lie in $D$. Hence, the Heegner point $R$ is not torsion.

Put $T=(12\omega^{i+2},-36)$
 which satisfies the relations
 \[T=\left\{ \begin{aligned} {[\omega^2]T+(0,12\sqrt{-3})},\\  [\omega]T+(0,-12\sqrt{-3}).\end{aligned}\right.\]
 Let $\alpha=2$ or  $1$ according to $p\equiv 4$ or $7\mod 9$ respectively. By Theorem \ref{co: Galois}, the point $Y=R-T$ belongs to $E_1(L_{(p)})^{\sigma_{\omega_3}=\omega^\alpha}$ which is identified with $E_{p}(K)$ under the isomorphism $(x,y)\mapsto ((\sqrt[3]{p})^2x,py)$ defined over $L_{(p)}$. Then the point $Y+\ov{Y}$ is identified with an element in $E_p(\BQ)$.
 \begin{prop}\label{divisibility}
 The point $Y$ is not divisible by $\sqrt{-3}$ in $E_1(L_{(p)})^{\sigma_{\omega_3}=\omega^\alpha}$. More precisely, there exists no point $X\in E_1(L_{(p)})^{\sigma_{\omega_3}=\omega^\alpha}$ and $S\in E_1(L_{(p)})^{\sigma_{\omega_3}=\omega^\alpha}_\tor$ such that 
 \[Y=\sqrt{-3}X+S.\]
 \end{prop}
\begin{proof}
Suppose
\begin{equation}\label{YX} 
Y=\sqrt{-3}X+S
\end{equation} 
for some $X\in E_1(L_{(p)})^{\sigma_{\omega_3}=\omega^\alpha}$ and $S\in E_1(L_{(p)})^{\sigma_{\omega_3}=\omega^\alpha}_\tor$. 
Since $H_{3p}/K$ is totally ramified at $\fP$ and $\ov{\fP}$, we have 

\[Y^{\sigma_{\omega_3}}\equiv Y,\ X^{\sigma_{\omega_3}}\equiv X,\ S^{\sigma_{\omega_3}}\equiv S\mod \fP  \textrm{  (resp. $\mod \ov{\fP}$)}.\]
From the formulae
\[Y^{\sigma_{\omega_3}}=[\omega^\alpha]Y,\quad X^{\sigma_{\omega_3}}=[\omega^\alpha]X,\quad S^{\sigma_{\omega_3}}=[\omega^\alpha]S,\]
we have 
\[[\sqrt{-3}]Y=[\sqrt{-3}]X=[\sqrt{-3}]S\equiv O\mod \fP  \textrm{  (resp. $\mod \ov{\fP}$)}.\]
By (\ref{YX}) and our choice of $T$, we have
\[S\equiv Y\equiv\left\{ \begin{aligned} {O,\quad\quad \quad \quad}\,&\mod \fP,\\ [1-\omega](12\omega^{i+1},-36), &\mod \ov{\fP},\end{aligned}\right.\]
which implies that $(x(S)\mod \fP,x(S)\mod \ov{\fP})$ doesn't lie in the subset $D$, and this contradicts with the fact that $S$ is a torsion point.
This proves that $Y$ is not divisible by $\sqrt{-3}$ in $E_1(L_{(p)})^{\sigma_{\omega_3}=\omega^\alpha}$.

\end{proof}
By the work of Dasgupta and Voight \cite{DV17}, we know that the free component of $E_p(K)$ has rank $1$ over $\CO_K$, and we have $$K\otimes_{\CO_K}E_p(K)\simeq K.$$ 
\begin{prop}\label{complex-conjugation}
As elements in $K\otimes_{\CO_K}E_p(K)$, we have 
\[\ov{Y}=\left(\omega^i\frac{\gamma}{\ov{\gamma}}\right)\otimes Y,\]
where $i=0,1,2$ and $\gamma \in \CO_K$ is a nonzero element satisfying $(\ov{\gamma},\gamma)=1$ and all primes of $\gamma$ are factors of rational primes which are split in $K$. 
\end{prop}
\begin{proof}
Recall if we set $\alpha=2$ or  $1$ according to $p\equiv 4$ or $7\mod 9$ respectively, then $E_1(L_{(p)})^{\sigma_{\omega_3}=\omega^\alpha}$ is identified with $E_p(K)$ under the real morphism 
\[\varphi: E_1\ra E_p,\quad  (x,y)\mapsto ((\sqrt[3]{p})^2x,py).\]
Since $R$ and $\ov{R}$ have the same height, as elements in 
 \[K\otimes _{\CO_K}E_1(L_{(p)})^{\sigma_{\omega_3}=\omega^\alpha},\]
 there exists $\beta\in K^\times $ such that
\[\N_{K/\BQ}(\beta)=1 \textrm{ and }\ov{R}=\beta\otimes R.\]
By Hilbert satz 90, we may assume $$\beta=u\cdot\frac{\gamma}{\bar{\gamma}},$$
where $u\in \CO_K^\times$, $\gamma \in \CO_K$ such that $(\ov{\gamma},\gamma)=1$ and all primes of $\gamma$ are factors of rational primes which are split in $K$. 

Then as usual points, there exists $S\in E_1(L_{(p)})^{\sigma_{\omega_3}=\omega^\alpha}_\tor=E_1(K)_\tor$ such that 
\[\ov{[\gamma]R}=[u\gamma]R+S.\]
Since $Y=\varphi(R-T)$, we have
\[\ov{[\gamma]Y}=[u\gamma]Y+S'.\]
for some $S'\in E_p(K)_\tor.$ 

It remains to prove that $u$ is a third root of unity. We have 
\begin{equation}\label{cc}
\ov{[\gamma]R}+[\gamma]R=[(u+1)\gamma]R+S.
\end{equation}
We claim that if $u$ is a primitive sixth root of unity or $-1$, then $[(u+1)\gamma]R$ has the same coordinates modulo both $\fP$ and $\ov{\fP}$. The case $u=-1$ is obvious. Suppose $u$ is a primitive six-th root of unity. Then $u+1=\sqrt{-3}\omega^j$ for some $j=1,2$. By Lemma \ref{lem:coord-R}, we have
\[
[u+1]R\equiv\left\{ \begin{aligned} {[\omega^{j+i+2}][\sqrt{-3}](12,-36) }&\mod \fP,\\ [\omega^{j+i+1}][\sqrt{-3}](12,-36)&\mod \ov{\fP}.\end{aligned}\right.
\]
Since $[\sqrt{-3}](12,-36)=(0,-12\sqrt{-3})$, we see that $[u+1]R$ has the same coordinates when modulo both $\fP$ and $\ov{\fP}$. Consequently, so does $[\gamma(1+u)]R$. Since $S$ is a $K$-point, the RHS of $(\ref{cc})$ has the same coordinates when modulo both $\fP$ and $\ov{\fP}$. 

On the other hand, we will show that the LHS of $(\ref{cc})$ has distinct coordinates when modulo $\fP$ and $\ov{\fP}$ respectively. To do this, it is enough to show that $\ov{[2\gamma]R}+[2\gamma]R$ has distinct coordinates when modulo $\fP$ and $\ov{\fP}$ respectively. Write $2\gamma=a+b\sqrt{-3}$ with $a,b\in\BZ$. Then
\[\ov{[2\gamma]R}+[2\gamma]R=[a](R+\ov{R})+[b\sqrt{-3}](R-\ov{R}).\]
By Lemma \ref{lem:coord-R}, $R+\ov{R}$ has distinct coordinates when modulo $\fP$ and $\ov{\fP}$ respectively and $R-\ov{R}\equiv [1-\omega^{i}]R$ for some $i=0,1,2$ when modulo $\fP$ and $\ov{\fP}$ respectively. 

Note that $R\mod \fP$ resp. $\ov{\fP}$ is of order $3$ in $E_1(\BF_p)$. Since $a\equiv 1,2\mod 3$, we know $[a](R+\ov{R})$ has distinct coordinates when modulo $\fP$ and $\ov{\fP}$ respectively. Since $\sqrt{-3}\mid (1-\omega^i)$, we conclude that 
\[ [b\sqrt{-3}](R-\ov{R})\equiv 0 \mod \fP \textrm{ resp. } \ov{\fP}.\]

So we conclude that that $\ov{[2\gamma]R}+[2\gamma]R$  and,  hence LHS of $(\ref{cc})$, have distinct coordinates when modulo $\fP$ and $\ov{\fP}$ respectively. Therefore, if $u$ is a primitive sixth root of unity or $-1$ we come to a contradiction and hence $u$ must be a third root of unity.
\end{proof}
\begin{remark}
If $p\equiv 7\mod 9$, it follows from Corollary \ref{co: Galois} that $$\ov{Y}\equiv [\omega^{\frac{p-7}{9}}]Y\mod \textrm{torsion}.$$ 
\end{remark}
\begin{thm}\label{divisibility2}
The point $Y+\ov{Y}\in E_p(\BQ)$ is not divisible by $3$.
\end{thm}
\begin{proof}
In the following all points are viewed as elements in $$K\otimes_{\CO_K}E_p(K)\simeq K.$$ 
By Proposition \ref{complex-conjugation}, there exists $\gamma \in \CO_K$ satisfying $(\ov{\gamma},\gamma)=1$ and all primes of $\gamma$ are factors of rational primes which are split in $K$ such that
\[Y+\ov{Y}=\left(1+u\frac{\gamma}{\ov{\gamma}}\right)\otimes Y.\]
Consider $\gamma+\ov{\gamma}$ as an element in $\CO_{K,3}=\BZ_3[\sqrt{-3}]$.  Suppose $\gamma=a+b\sqrt{-3}$ with $a,b\in \BZ_3$. Since  $\gamma$ is a $3$-adic unit, we have $a\in \BZ_3^\times$. Then 
\[1+u\frac{\gamma}{\ov{\gamma}}=1+\frac{\gamma}{\ov{\gamma}}+(u-1)\frac{\gamma}{\ov{\gamma}}=\frac{2a}{\gamma}+(u-1)\frac{\gamma}{\ov{\gamma}}.\]
is a $3$-adic unit since $\sqrt{-3}\mid (u-1)$. Then the 3-nondivisibility follows from Proposition \ref{divisibility}.
\end{proof}

\section{The explicit Gross-Zagier formulae}\label{E-G-Z}
\subsection{Test vectors and the explicit Gross-Zagier formulae}
Let $\pi$ be the automorphic representation of $\GL_2(\BA)$ corresponding to ${E_9}_{/\BQ}$. Then $\pi$ is only ramified at $3$ with conductor $3^5$. For $n\in \BQ^\times$, let $\chi_n: \Gal(K^{\ab}/K)\rightarrow\BC^\times$
be the cubic character given by $\chi_n(\sigma)=(\sqrt[3]{n})^{\sigma-1}$. Define
\[L(s,E_9,\chi_n):=L(s-1/2,\pi_K\otimes \chi_n),\quad \epsilon(E_9,\chi_n):=\epsilon(1/2,\pi_K\otimes \chi_n),\]
where $\pi_K$ is the base change of $\pi$ to $\GL_2(\BA_K)$.

Let $p\equiv 4,7\mod 9$ be a prime number, and put $\chi=\chi_{3p}$. From the Artin formalism, we have
$$L(s,E_9,\chi)=L(s,E_p)L(s,E_{3p^2}).$$
By \cite{Liverance}, we have  the epsilon factors $\epsilon(E_{3p^2})=+1$ and $\epsilon(E_p)=-1$, and hence the epsilon factor 
\[\epsilon(E_9,\chi)=\epsilon(E_{p})\epsilon(E_{3p^2})=-1.\] For a quaternion algebra $\BB_{\BA}$, we define its ramification index $\epsilon(\BB_v)=+1$ for any place $v$ of $\BQ$ if the local component $\BB_v$ is split and $\epsilon(\BB_v)=-1$ otherwise.
\begin{prop}\label{Tun-Saito}
The incoherent quaternion algebra $\BB$ over $\BA$, which satisfies
$$\epsilon(1/2,\pi_v,\chi_v)=\chi_{v}(-1)\epsilon_v(\BB)$$
for all places $v$ of $\BQ$, is only ramified at the infinity place.
\end{prop}
\begin{proof} Since $\pi$ is unramified at finite places $v\nmid 3$, $\chi$ is unramified at finite places $v\nmid 3p$ and $p$ is split in $K$, by \cite[Proposition 6.3]{Gross88} we get $\epsilon(1/2,\pi_v,\chi_v)=+1$ for all  finite $v\neq 3$. Again by \cite[Proposition 6.5]{Gross88}, we also know that $\epsilon(1/2,\pi_\infty,\chi_\infty)=-1$. Since $\epsilon(1/2,\pi,\chi)=-1$, we see that $\epsilon(1/2,\pi_3,\chi_3)=+1$. Since $\chi$ is a cubic character, $\chi_v(-1)=1$ for any $v$. Hence $\BB$ is only ramified at the infinity place.
\end{proof}
Let $\BB_f^\times=\GL_2(\BA_f)$ be the finite part of $\BB^\times $. For any open compact subgroup $U\subset \BB_f^\times$,  the Shimura curve $X_U$  associated to $\BB$ of level $U$ is the usual modular curve with complex uniformization
$$X_U(\BC)=\GL_2(\BQ)^+\backslash\left(\CH \bigsqcup \BP^1(\BQ)\right)\times \GL_2(\BA_f)/U.$$
Let $$\pi_{E_9}=\varinjlim_{U}\Hom^0_{\xi_U}(X_U,E_9),$$
where $\Hom^0_{\xi_{U}}(X_U,E_9)$ denotes the morphisms in $\Hom_\BQ(X_U,E_9)\otimes_\BZ\BQ$ using the Hodge class $\xi_U$ as a base point. Then $\pi_{E_9}$ is an automorphic representation of $\BB^\times$ over $\BQ$. Let $\pi$ be the Jacquet-Langlands correspondence of $\pi_{E_9}\otimes_\BQ\BC$ on $\GL_2(\BA)$.  By Proposition \ref{Tun-Saito} and a theorem of Tunnell-Saito \cite[Theorem 1.4.1]{YZZ}, the space
\[\Hom_{\BA_K^\times}(\pi_{E_9}\otimes \chi,\BC)\otimes \Hom_{\BA_K^\times}(\pi_{E_9}\otimes \chi^{-1},\BC)\]
is one-dimensional with a  generator $\beta=\otimes \beta_v$ where, for each place $v$ of $\BQ$, the bilinear form
\[\beta_v:\pi_{E_9,v}\otimes \pi_{E_{9,v}}\lra \BC\]
is given by
\begin{equation}\label{Eq:beta}
 \beta_v(\varphi_1,\varphi_2)=\int_{\BQ_v^\times\backslash K_v^\times}(\pi_{E_9,v}(t)\varphi_1,\varphi_2)\chi_v(t)dt,\quad \varphi_1,\varphi_2\in \pi_{E_9,v}.
\end{equation}
Here $(\cdot,\cdot)_v$ is a $\BB_v^\times$-invariant pairing on $\pi_{E_9,v}\times \pi_{E_9,v}$ and $dt$ is a Haar measure on $K_v^\times/\BQ_v^\times$.
For later application of the explicit Gross-Zagier formula in \cite{CST14}, if $(\varphi_{1},\varphi_{2})_v\neq 0$, we also define the normalized toric integral
\[\beta_v^0(\varphi_1,\varphi_2)=\int_{\BQ_v^\times\backslash K_v^\times}\frac{(\pi(t)\varphi_1,\varphi_2)_v\chi_v(t)}{(\varphi_1,\varphi_2)_v}dt.\]
 For more details we refer to \cite[Section 1.4]{YZZ} and \cite[Section 3]{CST14}.


The elliptic curve $E_{9}$ has conductor $3^5$. Let $f:X_0(3^5)\ra E_9$ be a nontrivial modular parametrization which sends the infinity cusp $[\infty]$ to the zero element $O$. Explicitly, we may take $f$ to be the quotient map $X_0(3^5)\ra X_\Gamma=E_9$ as in Proposition \ref{modular-curve}. Let $$\CR=\matrixx{\widehat{\BZ}}{\widehat{\BZ}}{3^5\cdot\widehat{\BZ}}{\widehat{\BZ}} \subset \BB_f(\wh{\BZ})=\M_2(\wh{\BZ})$$
be the Eichler order of discriminant $3^5$.
Then $U_0(3^5)=\CR^\times$, and by the newform theory \cite{Casselman1973}, the invariant subspace $\pi_{E_9}^{\CR^\times}$ has dimension $1$ and is generated by $f$.

\begin{prop}
The modular parametrization $f:X_0(3^5)\ra E_9$ is a test vector for the pair $(\pi_{E_9},\chi)$, i.e. $\beta(f,f)\neq 0$.
\end{prop}
\begin{proof}
Let $\CR'$ be the admissible order for the pair $(\pi_{E_9},\chi)$ in the sense of \cite[Definition 1.3]{CST14}. Since $\CR'$ and $\CR$ only differs at $3$. It suffices to verify that $\beta_3(f_3,f_3)\neq 0$, which is given by  Corollary \ref{ration}.
\end{proof}

Let $\omega_{E_n}$ be the invariant differential on the minimal model of $E_n$.  Define the minimal real period $\Omega_n$ of $E_n$ by
$$\Omega_{n}=\int_{E_n(\BR)}|\omega_{E_n}|.$$
By \cite[Formula (9)]{ZK}, we have
\begin{equation}\label{equation4}
\Omega_{p}\Omega_{3p^2}=(3p)^{-1}\Omega_9^2
\end{equation}
 Using Sagemath we compute that $\{\Omega_{9},\Omega_{9}\cdot(\frac{1}{2}+\frac{\sqrt{-3}}{2})\}$ is a $\BZ$-basis of the period lattice $L$ of the minimal model of $E_9$. So

\begin{equation}\label{period}
\sqrt{3}\Omega^2_{9}=2\int_{\BC/L}dxdy=\int_{E(\BC)}|\omega_{E_9}\wedge\overline{\omega}_{E_9}|=\frac{1}{6}\cdot 8\pi^2(\phi,\phi)_{\Gamma_0(3^5)},
\end{equation}
where $\phi$ is the newform of level $3^5$ and weight $2$ associated to $E_9$, and  $(\phi,\phi)_{\Gamma_0(3^5)}$ is the Petersson norm of $\phi$ defined by
$$(\phi,\phi)_{\Gamma_0(3^5)}=\int\int_{X_0(3^5)}|\phi(z)|^2dxdy,\ \ \ \ z=x+iy.$$
Recall $\tau=(2p\omega-9)/(9p\omega-36)\in \CH$ and let $P_1=[\tau,1]_{U_0(3^5)}$ be the CM point on $X_0(3^5)(H_{9p})$ and note that $f(P_1)=P_0$. Define the Heegner point
 \[R_1=\Tr_{H_{9p}/L_{(3,p)}} P_0\in E_9(L_{(3,p)}).\]
\begin{thm} \label{thm:GZ}
For primes $p\equiv 4,7 \mod 9$,  we  have the following explicit formula of Heegner points:
\[\frac{L'(1,E_{p})L(1,E_{3p^2})}{\Omega_{p}\Omega_{3p^2}}=2^{\alpha}\cdot 9 \cdot \wh{h}_\BQ(R_1),\]
where $\alpha=0$ if $p\equiv 4\mod 9$ and $\alpha=-1$ if $p\equiv 7\mod 9$.

\end{thm}
\begin{proof}
We note that the conductor of $\pi$ is $3^5$ and the conductor of $\chi$ is $9p$. Let $\CR'$ be the admissible order  for the pair $(\pi_{E_9},\chi)$ and let $f'\neq 0$ be a test vector in $V(\pi_{E_9},\chi)$ which is defined in \cite[Definition 1.4]{CST14}. The newform $f$ only differs from $f'$ at the local place $3$. Define the Heegner cycle
$$P^0_{\chi}(f)=\frac{\sharp \Pic(\CO_p)}{\Vol(\widehat{K}^\times/K^\times\widehat{\BQ}^\times,dt)}
\int_{K^\times\widehat{\BQ}^\times\backslash \widehat{K}^\times}f(P_1)^{\sigma_t}\chi(t)dt,$$
 and define $P^0_{\chi^{-1}}(f)$ similarly as in \cite[Theorem 1.6]{CST14}.
According to  Corollary \ref{ration},
\[\frac{\beta_3^0(f'_3,f'_3)}{\beta_3^0(f_3,f_3)}=2^{\alpha+2},\]
where $\alpha=0$ if $p\equiv 4\mod 9$ and $\alpha=-1$ if $p\equiv 7\mod 9$
. By \cite[Theorem 1.6]{CST14}, we have
$${L^{(3)}}'(1,E_9,\chi)=2^{\alpha+1}\cdot\frac{(8\pi^2)\cdot(\phi,\phi)_{\Gamma_0(3^5)}}{\sqrt{3}p\cdot(f,f)_{\CR'}}\cdot\left\langle P^0_{\chi}(f), P^0_{\chi^{-1}}(f)\right\rangle_{K,K},$$
where ${L}^{(3)}$ denotes the partial $L$-function with the $3$-adic local factor removed, $(\cdot,\cdot)_{\CR'}$ is the the pairing on $\pi_{E_9}\times \pi_{{E_9}^\vee}$ defined as in \cite[page 789]{CST14}, and $\langle\cdot,\cdot\rangle_{K,K}$ is a pairing from $E_9(\ov{K})_\BQ\times_K E_9(\ov{K})_\BQ$ to $\BC$  such that $\langle\cdot,\cdot\rangle_{K}=\rm{Tr}_{\BC/\BR}\langle\cdot,\cdot\rangle_{K,K}$ is the Neron-Tate height over the base field $K$, see \cite[page 790]{CST14}.  The local representation $\pi_{K,3}\otimes \chi_3$ is the principal series induced from the pair $(\Theta_3\chi_3,\ov{\Theta_3}\chi_3)$, where $\Theta$ is the Hecke character over $K$ associated to $E_9$ via the CM theory. By Lemma \ref{thetavalue} and \ref{chi}, both the characters $\Theta_3\chi_3,\ov{\Theta_3}\chi_3$ are ramified. Hence the $3$-adic local $L$-factor of $L(s,E_9,\chi)$ is trivial and we have
\[{L'}(1,E_9,\chi)={L^{(3)}}'(1,E_9,\chi).\]
In our case, by \cite[Lemma 2.2, Lemma 3.5]{CST14},
$$(f,f)_{\CR'}=\frac{\Vol(X_{\CR^{' \times}})}{\Vol(X_{\CR^\times})}\deg f=6\cdot \frac{\Vol(\CR^\times)}{\Vol(\CR^{'\times})}=4.$$
So, we get
\begin{equation}
L'(1,E_9,\chi)=2^{\alpha-1}\frac{(8\pi^2)\cdot(\phi,\phi)_{\Gamma_0(3^5)}}{\sqrt{3p^2}}\cdot\left\langle P^0_{\chi}(f), P^0_{\chi^{-1}}(f)\right\rangle_{K,K}.
\end{equation}

On the other hand
$$P^0_{\chi}(f)=\frac{\# \Pic(\CO_p)}{\# \Pic(\CO_{9p})}\sum_{t\in\Pic(\CO_{9p})}f(P_1)^{\sigma_t}\chi(t).$$
Since
$$\frac{\#\Pic(\CO_p)}{\#\Pic(\CO_{9p})}=[K^\times \wh{\CO}_p^\times:K^\times \wh{\CO}_{9p}^\times]^{-1}=\frac{1}{9},$$
we have
$$P^0_{\chi}(f)=\frac{1}{9}\sum_{t\in\Pic(\CO_{9p})}f(P_1)^{\sigma_t}\chi(t).$$
If we put
\[R_2=\sum_{\sigma\in \Gal(H_{9p/L(3p)})} f(P_1)^\sigma\chi(\sigma)=3R_1\in E_9(L_{(3p)}),\]
then
\begin{align}\langle P^0_{\chi}(f),P^0_{\chi^{-1}}(f)\rangle_{K,K}&=\frac{1}{9^2}
\langle\sum_{\sigma\in\Gal(L_{(3p)}/K)}R_2^{\sigma}\chi(\sigma),\sum_{\sigma\in\Gal(L_{(3p)}/K)}
R_2^{\sigma}\chi^{-1}(\sigma)\rangle_{K,K}\\
&=\frac{1}{27}\langle R_2,\sum_{\sigma\in\Gal(L_{(3p)}/K)}R_2^{\sigma}\chi^{-1}(\sigma)\rangle_{K,K}\notag\\
&=\frac{1}{27}(\langle R_2,R_2\rangle_{K,K}+\chi^{-1}(\sigma')\langle R_2,R_2^{\sigma'}\rangle_{K,K}+\chi^{-1}(\sigma'^2)\langle R_2,R_2^{\sigma'^2}\rangle_{K,K})\\
&=\frac{1}{27}\left(\langle R_2,R_2\rangle_{K,K}-\left\langle R_2,R_2^{\sigma'}\right\rangle_{K,K}\right),\notag
\end{align}
where $\sigma'$ is a generator of $\Gal(L_{(3p)}/K)$. In the last equality we use the fact that $\langle R_2,R_2^{\sigma'}\rangle_{K,K}=\langle R_2,R_2^{\sigma'^2}\rangle_{K,K}$ since $\langle,\rangle_{K,K}$ is symmetric and Galois invariant. 
By Theorem \ref{thm:Galois} and Corollary \ref{co: Galois}, we can assume $R_2^{\sigma'}=[\omega]R_2$, then
$$\left\langle R_2,R_2^{\sigma'}\right\rangle_{K,K}=\frac{1}{2}\left(\widehat{h}_K([1+\omega]R_2)
-\widehat{h}_K([\omega]R_2)-\widehat{h}_K(R_2)\right).$$
Since $|1+\omega|=|\omega|=1$, by definition, $\widehat{h}_K([1+\omega]R_2)=\widehat{h}_K([\omega]R_2)=\widehat{h}_K(R_2)$. Then
$$
\left\langle R_2,R_2^{\sigma'}\right\rangle_{K,K}=-\frac{1}{2}\widehat{h}_K(R_2),$$
and hence
\begin{equation}\label{equation2}
\left\langle P^0_{\chi}(f),P^0_{\chi^{-1}}(f)\right\rangle_{K,K}=\frac{1}{18}\widehat{h}_K(R_2)=\frac{1}{9}\widehat{h}_\BQ(R_2)=\wh{h}_\BQ(R_1).
\end{equation}

Finally, combining  (\ref{equation4})-(\ref{equation2}), we get
\[\frac{L'(1,E_{p})L(1,E_{3p^2})}{\Omega_{p}\Omega_{3p^2}}=2^{\alpha}\cdot 9 \cdot \wh{h}_\BQ(R_1).\]
\end{proof}

Recall there is an isomorphism
\[\phi:E_9\lra E_1,\quad (x,y)\mapsto \left(\left(\sqrt[3]{3}\right)^2x, 3y\right),\]
and we have the following commutative diagram:
$$\xymatrix{E_9(H_{9p})^{\sigma_{1+3\omega_3}=\omega^2}\ar[d]^\phi\ar[rr]^{\Tr_{H_{9p}/L_{(3,p)}}}&&E_9(L_{(3,p)})^{\sigma_{1+3\omega_3}=\omega^2}\ar[d]^{\phi}\\
            E_1(H_{3p})\ar[rr]^{\Tr_{H_{3p}/L_{(p)}}}&&E_1(L_{(p)})}.$$
In paricular, we have $\phi(R_1)=R$, and hence the following
\begin{coro}\label{co:GZ}
For primes $p\equiv 4,7 \mod 9$,  we  have
\[\frac{L'(1,E_{p})L(1,E_{3p^2})}{\Omega_{p}\Omega_{3p^2}}=2^{\alpha}\cdot 9 \cdot \wh{h}_\BQ(R),\]
where $\alpha=0$ if $p\equiv 4\mod 9$ and $\alpha=-1$ if $p\equiv 7\mod 9$.
\end{coro}
\begin{proof}
This is immediate from Theorem \ref{thm:GZ}.
\end{proof}
Recall that Dasgupta and Voight \cite{DV17} proved that the Heegner point $R$ is not torsion. By the above Gross-Zagier formula and the work of Kolyvagin \cite{Kolyvagin1990}, we know that
\[\rk_\BZ E_p(\BQ)=\ord_{s=1}L(s,E_p)=1,\quad \rk_\BZ E_{3p^2}(\BQ)=\ord_{s=1}L(s,E_{3p^2})=0.\]

\section{The $3$-part of the Birch and Swinnerton-Dyer conjectures}\label{E-B-SD}
Let $F$ be a number field. Let $\phi:A\rightarrow A'$ be an isogeny of elliptic curves over $F$ of degree $m$ and $\phi'$ is its dual isogeny. The commutative diagram
$$\xymatrix{0\ar[r]&A[\phi]\ar[d]\ar[r]&A\ar[r]^{\phi}\ar@{=}[d]&A'\ar[r]\ar[d]^{\phi'}&0\\
            0\ar[r]&A[m]\ar[d]^\phi\ar[r]&A\ar[r]^{[m]}\ar[d]^{\phi}&A\ar[r]\ar@{=}[d]&0\\
            0\ar[r]&A'[\phi']\ar[r]&A'\ar[r]^{\phi'}&A\ar[r]&0}$$
induces the following commutative diagram
\[
\xymatrix@C=7pt{&0\ar[d]&&0\ar[d]&&& \\
0\ar[r]&A'[\phi'](F)/\phi A[m](F)\ar[d]^{\phi'}\ar[rr]^{=}&&A'[\phi'](F)/\phi A[m](F)\ar[d]^{\phi'}\ar[rr]&&0\ar[d]&\\
0\ar[r]&A'(F)/\phi A(F)\ar[d]^{\phi'}\ar[rr]&&\Sel_\phi(A/F)\ar[rr]\ar[d]&&\Sha
(A/F)[\phi]\ar[r]\ar[d]&0\\
0\ar[r]&A(F)/mA(F)\ar[d]\ar[rr]&&\Sel_m(A/F)\ar[rr]\ar[d]^{\phi}&&\Sha
(A/F)[m]\ar[r]\ar[d]^{\phi}&0\\
0\ar[r]&A(F)/\phi'A'(F)\ar[d]\ar[rr]&&\Sel_{\phi'}(A'/F)\ar[rr]\ar[d]&&\Sha(A'/F)[\phi']\ar[r]\ar[d]&0\\
&0\ar[rr]&&\Sel_{\phi'}(A'/F)/\phi \Sel_m(A/F)\ar[rr]^{\simeq}\ar[d]&&\Sha(A'/F)[\phi']/\phi\Sha(A/F)[m]\ar[r]\ar[d]&0\\
&&&0&&0&.
}
\]
From this diagram, we immediately have the following
\begin{lem} \label{Sel}
Let $A,A'$ and $\phi,\phi'$ be as above.
 \[
|\Sel_m(A/F)|=\frac{|\Sel_\phi(A/F)||\Sel_{\phi'}(A'/F)|}{|A'[\phi'](F)/\phi A[m](F)||\Sha(A'/F)[\phi']/\phi\Sha(A/k)[m]|}.
\]
\end{lem}

Let $n$ be a positive cube-free integer and $E'_n$ be the elliptic curve given by Weierstrass equation $y^2=x^3+(4n)^2$. Then there is an unique isogeny
$\phi_n: E_n\rightarrow E'_n$ of degree $3$ up to $[\pm 1]$ and denote $\phi_n'$ its dual isogeny.

\begin{prop}\label{sha}
Let $p\equiv 4,7\mod 9$ be primes such that $3\mod p$ is not a cubic residue. Then
$$\dim_{\BF_3}\Sel_3(E_p(\BQ))\leq 1,\quad \dim_{\BF_3}\Sel_3(E_{3p^2}(\BQ))=0.$$
\end{prop}
\begin{proof}
By \cite[Theorem 2.9]{Satge}, we know that
$$\Sel_{\phi_p}(E_p(\BQ))=\Sel_{\phi_p'}(E'_p(\BQ))=\BZ/3\BZ,$$
and
$$\Sel_{\phi_{3p^2}}(E_{3p^2}(\BQ))=\BZ/3\BZ,\quad \Sel_{\phi_{3p^2}'}(E'_{3p^2}(\BQ))=0.$$
Note that $E_p[3](\BQ)$ and $E_{3p^2}[3](\BQ)$ are trivial and $|E'_p[\phi'_p](\BQ)|=|E'_{3p^2}[\phi'_{3p^2}](\BQ)|=3$. By Lemma \ref{Sel}, the proposition follows.
\end{proof}

Now we are ready to give the proof of Theorem \ref{Main}.
\begin{proof}[Proof of Theorem 1.4]
By \cite[Table 1]{ZK}, we know that $c_p(E_p)=3$, $c_3(E_p)=1\ \rm{or}\ 2$ depending on $p$ congruent to $4$ or $7$ modulo $9$ respectively, and $c_\ell(E_p)=1$ for primes $\ell\neq 3,p$, while $c_p(E_{3p^2})=3$, $c_\ell(E_{3p^2})=1$ for primes $\ell\neq p$. 

Let $P$ be the generator of the free part of $E_p(\BQ)$. Then the BSD conjecture predicts that
\begin{equation*}
\frac{L'(1,E_p)}{\Omega_p}=2^m\cdot 3\cdot |\Sha(E_p)|\cdot\widehat{h}_\BQ(P),
\end{equation*}
where $m=0$ if $p\equiv 4\mod 9$ and $1$ if $p\equiv 7\mod 9$, and
\begin{equation*}
\frac{L(1,E_{3p^2})}{\Omega_{3p^2}}=3\cdot\left|\Sha(E_{3p^2})\right|.
\end{equation*}
Combining these two, we get
\begin{equation*}
\frac{L'(1,E_p)}{\Omega_p}\cdot\frac{L(1,E_{3p^2})}{\Omega_{3p^2}}=2^m\cdot 9\cdot|\Sha(E_p)|\cdot|\Sha(E_{3p^2})|\cdot\widehat{h}_\BQ(P).
\end{equation*}
By Theorem \ref{thm:GZ} and Corollary \ref{co:GZ}, we expect
\begin{equation}\label{bsd}
|\Sha(E_p)|\cdot|\Sha(E_{3p^2})|=2^i\frac{\wh{h}_\BQ(R)}{\wh{h}_\BQ(P)},
\end{equation}
where $i=0$ resp. $i=-2$ if $p\equiv 4\mod 9$ resp. $p\equiv 7\mod 9$. Note the RHS of $(\ref{bsd})$
is a nonzero rational number. 

By Proposition \ref{sha}, $E_p$ being rank 1, and the exact sequence
\[\xymatrix{0\ar[r]&E(\BQ)/3E(\BQ)\ar[r]&Sel_3(E(\BQ))\ar[r]&\Sha(E)[3]\ar[r]&0},\]
we know directly that
\begin{equation}
 |\Sha(E_p)[3^\infty]|=|\Sha(E_{3p^2})[3^\infty]|=1.
\end{equation}
In order to prove the $3$-part of (\ref{bsd}), it suffices to prove
$$\widehat{h}_\BQ(P)=u\widehat{h}_\BQ(R_1)=u\widehat{h}_\BQ(R)$$
for some $u\in \BZ_3^\times \cap \BQ$.  However, it follows from Proposition \ref{divisibility} and Theorem \ref{divisibility2} that
\[\wh{h}_\BQ(P)=w\wh{h}_\BQ(\ov{Y}+Y)=w\wh{h}_\BQ(\ov{R}+R)=u\wh{h}_\BQ(R)\]
for some $u,w\in\BZ_3^\times \cap \BQ$, where the last equality follows from Lemma \ref{complex-conjugation}. Indeed, we note that a similar formula for the complex conjugation is valid for $R$ in the proof of Lemma \ref{complex-conjugation}.
\end{proof}

\section{Waldspurger's period integral of minimal vectors}\label{minimal vector}
This whole section is purely local, so we omit subscript $v$ from all notations.
Let $\pi$ be a smooth irreducible representation of $\GL_2$ over a $p$-adic field $\BF$ with central character $w_\pi$, and $\chi$ be any character over a quadratic extension $\BE/\BF$ such that  $w_\pi \chi|_{\BF^\times}=1$. 
For any test vector $\varphi\in\pi$, denote the local Waldspurger's period integral against a character $\chi$ on $\BE^\times$ by
\begin{equation}
\WaldsI=\int\limits_{t\in \BF^\times\backslash \BE^\times} \Phi_{\varphi}(t)\chi(t)dt
\end{equation}
where $\Phi_{\varphi}=(\pi(t)\varphi,\varphi)$ is the matrix coefficient associated to $\varphi$. 

In this section we shall start with the general setting, review the compact induction theory and minimal vectors which arise naturally, and then recall the main results from \cite{HN18} which evaluates $\WaldsI$ when $\varphi$ is a  minimal vector.

\subsection{Preliminaries}\label{prelimi}
For a real number $a$, let $\lfloor a\rfloor\leq a$ be the largest possible integer, and $\lceil a\rceil \geq a$ be the smallest possible integer.

Let $\F$ be a $p$-adic field with \New{residue field of order} $q$\yhb{residue field of order $q$ better?}, uniformizer $\varpi=\varpi_\F$, ring of integers $\CO_\F$ and $p$-adic valuation $v_\F$. Let $\psi$ be an additive character of $\F$.\yhb{Maybe it is better to move the sentence before the definition of $\psi_\BL$ below.}\yh{It's according to the field. This is a minor issue anyway.} Assume that $2\nmid q$. For $n\geq 1$, let $U_\F(n)=1+\varpi_\F^n \CO_\F$.
Let $\pi$ be a supercuspidal representation over $\F$ with central character $w_\pi=1$. 

Let $\L$  be a quadratic extension over $\F$. Let $e_\L=e(\L/\F)$ be the ramification index and $v_\L$ be the valuation on $\L$. Let $\varpi_\L$ be a uniformizer for $\L$. When $\L$ is unramified we shall identify $\varpi_\L$ with $\varpi_\F$. Otherwise we suppose that $\varpi_\L^2=\varpi_\F$. Let $x\mapsto \overline{x}$ to be the unique nontrivial involution of $\L/\F$.
 Let $\psi_\L=\psi\circ \Tr_{\L/\F}$. One can make similar definitions for a possibly different quadratic extension $\E$. Note that we shall assume that $\varpi_\E^2=\xi\varpi_\F$ for $\xi\in \CO_\F^\times-(\CO_\BF^\times)^2$ if $\E$ and $\L$ are both ramified and distinct.

For $\chi$ a multiplicative character on $\F^\times$, let $c(\chi)$ be  the smallest integer such that $\chi$ is trivial on $1+\varpi_\F^{c(\chi)} \CO_\F$. Similarly $c(\psi)$ is the smallest integer such that  $\psi$ is trivial on $\varpi_\F^{c(\psi)}\CO_\F$. We choose $\psi$  to be unramified, or equivalently,   $c(\psi)=0$. Then $c(\psi_\L)=-e_\L+1$.
Let $c(\pi)$ be the power of the conductor of $\pi$.

When $\chi$ is a character over a quadratic extension, denote $\overline{\chi}(x)=\chi(\overline{x})$.
\Cor{\begin{lem}
Suppose that $p$ is large enough.
For a multiplicative character $\nu$ over $\F$ with $c(\nu)\geq 2$, there exists $\alpha_\nu\in \F^\times$ with $v_\F(\alpha_\nu)=-c(\nu)+c(\psi)$ such that for $u\in \varpi_\F \CO_\F$
\begin{equation}
\nu(1+u)=\psi(\alpha_\nu \log(1+u)) \text{\yhb{what is $\psi_\BF$?}}\yh{\text{Changed to }\psi}
\end{equation}
where $\log(1+u)$ is the standard Taylor expansion for logarithm 
\begin{equation}
\log(1+u)=u-\frac{u^2}{2}+\frac{u^3}{3}+\cdots.
\end{equation}
In particular when $p>2$ we have that 
\begin{equation}
\nu(1+u)=\psi(\alpha_\nu (u-\frac{u^2}{2})) \end{equation}
for any $u\in \varpi_\F^{\lfloor c(\nu)/2\rfloor} \CO_\F$.
\end{lem}}
\New{\begin{lem}\label{Lem:DualLiealgForChar}
For a multiplicative character $\nu$ over $\F$ with $c(\nu)\geq 2$, there exists $\alpha_\nu\in \F^\times$ with $v_\F(\alpha_\nu)=-c(\nu)+c(\psi)$ such that
\begin{equation}\label{eq:alphatheta}
\nu(1+u)=\psi(\alpha_\nu u) \end{equation}
for any $u\in \varpi_\F^{\lceil c(\nu)/2\rceil} \CO_\F$. $\alpha_\nu$ is determined $\mod \varpi_\F^{-\lceil c(\nu)/2\rceil+c(\psi)} \CO_\F$.
\end{lem}
One can easily check this lemma by using that $\nu(1+u)$ becomes an additive character in $u$ for $u\in \varpi_\F^{\lceil c(\nu)/2\rceil} \CO_\F$.
}

We shall also need some basic results for compact induction theory and matrix coefficient.
In general let $G$ be a unimodular locally profinite group with center $Z$. Let $H\subset G$ be an open and closed subgroup containing $Z$ with $H/Z$ compact. Let $\rho$ be an irreducible smooth representation of $H$ with unitary central character and 
\New{$$\pi=c-\Ind_H^G(\rho)=\{f:G\rightarrow \rho| f(hg)=\rho(h)f(g)\forall h\in H, \text{\ $f$ is compactly supported}\}.$$}
\yhb{explain the notation, Maybe better to give the definition of compact induction?}  By the assumption on $H/Z$, $\rho$ is automatically unitarisable, and we shall denote the unitary pairing on $\rho$ by $<\cdot,\cdot>_{\rho}$. Then one can define the unitary pairing on $\pi$ by
\begin{equation}
<\phi,\psi>=\sum\limits_{x\in H\backslash G}<\phi(x),\psi(x)>_{\rho}.
\end{equation}
If we let $y\in H\backslash G$ and $\{v_i\}$ be a basis for $\rho$, \yhb{orthnormal basis?}\yh{Not necessary}the elements
\begin{equation}\label{Eq:CompactIndBasis}
 f_{y,v_i}(g)=\begin{cases}
\rho(h)v_i,&\text{\ if \ }g=hy\in Hy;\\
0,&\text{\ otherwise.}
\end{cases}
\end{equation} 
form a basis for $\pi$.
\begin{lem}\label{lem3.3:matrixcoeffInduction}
For $y,z\in H\backslash G$,
\begin{equation}
<\pi(g)f_{y,v_i},f_{z,v_j}>=\begin{cases}
<\rho(h)v_i,v_j>_{\rho}, &\text{\ if\ }g=z^{-1}hy\in z^{-1}Hy;\\
0,&\text{\ otherwise}.
\end{cases}
\end{equation}
\end{lem}
\yhb{Please give the reference of Lemma 6.1 and Lemma 6.2, or skecth the proof}

\subsection{Compact induction theory for supercuspidal representations and minimal vectors}
We shall review the compact induction theory for supercuspidal representations on $\GL_2$. For more details, see \cite{BushnellHenniart:06a}.

We shall fix the embeddings and work out everything explicitly. For $v_\F(D')=0,1$,\yhb{make it clear $v$ is the valuation over $\BF$} we shall refer to the following embedding of a quadratic field extension $\BL=\BF(\sqrt{D'})$ as a standard embedding.
\begin{equation}\label{Eq:standardembeddingL}
   x+y\sqrt{D'}\mapsto \matrixx{x}{y}{yD'}{x}.
\end{equation}

Supercuspidal representations are parametrised via compact induction by characters $\theta$ over some quadratic extension $\L$. More specifically we have the following quick guide. \begin{enumerate}
\item[Case 1.]$c(\pi)=2n+1$ corresponds to $e_\BL=2$ and $c(\theta)=2n$ .
\item[Case 2.] $c(\pi)=4n$ corresponds to $e_\BL=1$ and $c(\theta)=2n$.
\item[Case 3.] $c(\pi)=4n+2$ corresponds to $e_\BL=1$ and $c(\theta)=2n+1$ .
\end{enumerate}

\begin{defn}
For $e_{\BL}=1,2$, let $$\fA_{e_{\BL}}=\begin{cases}
M_2(\CO_{\BF}), \text{\ if }e_\BL=1,\\
\matrixx{ \CO_\BF}{\CO_\BF}{\varpi \CO_\BF}{ \CO_\BF},\text{\ otherwise}.
\end{cases}
$$
Its Jacobson radical is given by
$$\CB_{e_\BL}=\begin{cases}
\varpi M_2(\CO_{\BF}), \text{\ if }e_\BL=1,\\
\matrixx{\varpi \CO_\BF}{\CO_\BF}{\varpi \CO_\BF}{\varpi \CO_\BF},\text{\ otherwise}.
\end{cases}$$
\end{defn}

Define  filtrations of compact open subgroups as follows:
\begin{equation}
K_{\fA_{e_\BL}}(n)=1+\CB_{e_\BL}^n,\text{\yhb{what is $\CB$?}\yh{Added subscript $e_\BL$}} U_{\BL}(n)=1+\varpi_{\BL}^n\CO_{\BL}.
\end{equation}

Note that each $K_{\fA_{e_\BL}}(n)$ is normalised by $\BL^\times$ which is embedded as in \eqref{Eq:standardembeddingL}.

Denote 
$J=\BL^\times K_{\fA_{e_\BL}}(\lfloor c(\theta)/2\rfloor)$, $J^1=U_\BL(1)K_{\fA_{e_\BL}}(\lfloor c(\theta)/2\rfloor)$, $H^1=U_\BL(1)K_{\fA_{e_\BL}}(\lceil c(\theta)/2\rceil)$. Then $\theta$ on $\BL^\times$ can be extended to be a character $\tilde{\theta}$ on $H^1$ by
\begin{equation}\label{Eq:thetatilde}
\tilde{\theta}(l (1+x))=\theta(l)\psi\circ\Tr (\alpha_\theta x),
\end{equation}
where $l\in \BL^\times$, $1+x\in K_{\fA_{e_\BL}}(\lceil c(\theta)/2\rceil)$ and $\alpha_\theta\in \BL^\times\subset M_2(\BF)$ is associated to $\theta$ by Lemma \ref{Lem:DualLiealgForChar} under the fixed embedding.

When $c(\theta)$ is even, $H^1=J^1$ and $\tilde{\theta}$ can be further extended to $J$ by the same formula. In this case denote $\Lambda=\tilde{\theta}$  and $\pi=c-\Ind_J^G\Lambda$.

When $c(\theta)$ is odd, $J^1/H^1$ is a two dimensional vector space over the residue field. This case only occurs when $c(\pi)=4n+2$ as listed above.
Then there exists a $q-$dimensional representation $\Lambda$ of $J$ such that
$\Lambda|_{H^1}$ is a multiple of $\tilde{\theta}$, and
$\Lambda|_{\BL^\times}=\oplus \theta\nu$ where $c(\nu)=1$ and $\nu|_{\BF^\times}=1$. More specifically, let $B^1$ be any intermediate group between $J^1$ and $H^1$ such that $B^1/H^1$ gives a polarisation of $J^1/H^1$ under the  pairing given by \Cor{$\psi_{\alpha_\theta}$.}
\New{$$(1+x,1+y)\mapsto \psi\circ \Tr (\alpha_\theta [x, y]).$$}
\yhb{what is $\psi_{\alpha_\theta}$?} Then $\tilde{\theta}$ can be extended to $B^1$ by the same formula \eqref{Eq:thetatilde} and  $\Lambda|_{J^1}=\Ind_{B^1}^{J^1}\tilde{\theta}$.

 In the case $J^1=H^1$, we take $B^1=J$ for uniformity. In either cases, we have $w_\pi=\theta|_{\BF^\times}$.

\begin{defn}
There exists a unique element $\varphi_0\in \pi$ such that $B^1$ acts on it by $\tilde{\theta}$. (Type 1 minimal vector in the notation of \cite{HN18}.)

We also call any single translate $\pi(g)\varphi_0$ a minimal vector, as the conjugated group $gB^1g^{-1}$ acts on it by the conjugated character $\tilde{\theta}^g$. They form an orthonormal basis for $\pi$.
\end{defn}
\begin{coro}\label{Cor:MCofGeneralMinimalVec}
Let $\Phi_{\varphi_0}$ be the matrix coefficient associated to a minimal vector $\varphi_0$ as above. Then $\Phi_{\varphi_0}$ is supported on $J$, and
\begin{equation}
\Phi_{\varphi_0}(bx)=\Phi_{\varphi_0}(xb)=\tilde{\theta}(b)\Phi_{\varphi_0}(x)
\end{equation}
for any $b\in B^1$. Furthermore when $\dim \Lambda\neq 1$, $\Phi_{\varphi_0}|_{J^1}$ is supported only on $B^1$.
\end{coro}
Note that $\varphi_0$ is basically $f_{1,v_i}$ as in \eqref{Eq:CompactIndBasis} for the coset representative $1\in J\backslash G$. The corollary \New{follows immediately}\yhb{immediate or immediately?} from the definition of $\varphi_0$ and Lemma \ref{lem3.3:matrixcoeffInduction}.

\subsection{Local Langlands correspondence and compact induction}\label{SubSec:Langlands-compactInd}
Here we describe the relation between the compact induction parametrisation and the local Langlands correspondence. See \cite{BushnellHenniart:06a} Section 34 for more details.

For a field extension $\L/\F$ and an additive character $\psi$ over $\F$, let $\lambda_{\L/\F}(\psi)$ be the Langlands $\lambda-$function in \cite{Langlands}. When $\L/\F$ is a quadratic field extension, let $\eta_{\L/\F}$ be the associated quadratic character. By \cite[Lemma 5.1]{Langlands}, we have for $\psi_\beta(x)=\psi(\beta x)$,
\begin{equation}\label{Eq:LanglandsLambdaFun}
\lambda_{\L/\F}(\psi_\beta)=\eta_{\L/\F}(\beta)\lambda_{\L/\F}(\psi).
\end{equation}
\begin{defn}\label{Defn:DeltaTheta}

\begin{enumerate}
\item If $\BL$ is inert, define $\Delta_\theta$ to be the unique unramified character of $\BL^\times$ of order 2.\yhb{what does the character look like?}\yh{Take -1 on uniformizer}

\item If $\BL$ is ramified and $\theta$ is a character over $\BL$ with $c(\theta)>0$ even, associate $\alpha_\theta$ to $\theta$ as in Lemma \ref{Lem:DualLiealgForChar}.
Then define $\Delta_\theta$ to be the unique level 1 character of $\BL^\times$ such that
\begin{align}
\Delta_\theta|_{\BF^\times}=\eta_{\BL/\BF}, \Delta_\theta(\varpi_\BL)=\eta_{\BL/\BF}(\varpi_\BL^{c(\theta)-1}\alpha_\theta)\lambda_{\BL/\BF}^{c(\theta)-1}(\psi).
\end{align}
\end{enumerate}

\end{defn}
Note that in  \cite{BushnellHenniart:06a}   $\psi$ is chosen to be level 1. We have adapted the formula there to our choice of $\psi$ using \eqref{Eq:LanglandsLambdaFun}.
The definition is also independent of the choice of $\varpi_\BL$.    
\begin{thm}\label{LLC}
If $\pi$ is associated by compact induction to a character $\theta$ over a quadratic extension $\BL$, then its associated Deligne-Weil representation by local Langlands correspondence is $\sigma=\Ind_{\BL}^{\BF} (\Theta)$, where $\Theta=\theta\Delta_\theta^{-1}$, or equivalently $\theta=\Theta\Delta_{\Theta}$.
\end{thm}
Note here that $\Theta$ and $\theta$ always differ by a  level $\leq 1$ character, so $\alpha_\Theta$ \New{can be chosen to be the same as} $\alpha_\theta$ in Lemma \ref{Lem:DualLiealgForChar}\yhb{modulo what?}and $\Delta_\Theta=\Delta_\theta$.

\subsection{Using minimal vectors for Waldspurger's period integral }\label{Sec:WaldsForMinimalVec}
\New{Now we review the local Waldspurger's period integral for minimal vectors. For details and proofs, see \cite{HN18} or the appendix.}

For simplicity we pick $$D'=\frac{1}{\alpha_\theta^2\varpi_\L^{2c(\theta)}},$$
identify $\frac{1}{\alpha_\theta\varpi_\L^{c(\theta)}}$ with $\sqrt{D'}$ and $\L$ with $\F(\sqrt{D'})$, and use the standard embedding \eqref{Eq:standardembeddingL}.

By this choice, we have $v_\F(D')=0$ if $e_\L=1$ and $v_\F(D')=1$ if $e_\L=2$.
\begin{equation}\label{eq2.1:specialembedding}
\alpha_\theta= \frac{1}{\varpi_\L^{c(\theta)}}\frac{1}{\sqrt{D'}} \mapsto \frac{1}{\varpi^{c(\theta)/e_\L}}\zxz{0}{\frac{1}{D'}}{1}{0}.
\end{equation}
Such choice is not essential. A different choice will result in slightly different equations, for example, \eqref{eq:new-necessary-same-ramified} and \eqref{eq:3.4:IntertwiningtoWhittaker}, but the final results are similar.

Note that when $e_\L=2$, $c(\theta)$ must be even when $\theta|_{\F^*}=1$. 
 Assume that $\BE=\BF(\sqrt{D})$ for $v_\F(D)=0,1$ is also embedded as
\begin{equation}\label{Eq:standardembedding}
   x+y\sqrt{D}\mapsto \matrixx{x}{y}{yD}{x}. 
\end{equation}

In \cite{HN18}, test vectors of form $\pi(g)\varphi_0$ were used to study $\WaldsI$ for general combinations of $\pi$, $\BE$ and $\chi$. 
For the purpose of this paper we shall only review the case when $\BL\simeq\BE$ are ramified,  $w_\pi=1$, $c(\pi)=2n+1$ is odd and $c(\pi_{\chi})\leq c(\pi)$. 
We shall normalise the Haar measure on $\BE^\times$ so that $\Vol( \CO_\BF^\times\backslash \CO_\BE^\times)=1$. Then $\Vol(\BF^\times\backslash \BE^\times)=2$.

In this case we use test vectors of form  $$\varphi=\pi\lrb{\matrixx{1}{u}{0}{1}\matrixx{v}{0}{0}{1}}\varphi_0$$
for some $u\in \CO_\BF,v\in \CO_\BF^\times$. Recall that when $e_\BL=2$, 
$$K_{\fA_{e_\BL}}(n)=1+\matrixx{\varpi^{ \lceil n/2\rceil}\CO_\BF}{\varpi^{ \lfloor n/2\rfloor}\CO_\BF}{\varpi^{ \lfloor n/2\rfloor+1}\CO_\BF}{\varpi^{ \lceil n/2\rceil}\CO_\BF},$$
and $J=\BL^\times K_{\fA_{e_\BL}}(n)$ acts on $\varphi_0$ by a character, so we can assume that $v\in (\CO_\BF/ \varpi^{\lceil n/2\rceil}\CO_\BF)^\times$.

There are two situations depending on $\min \{c(\theta\chi), c(\theta\overline{\chi})
\}$. Note that for the embedding of $\BE$ fixed above, we have
\begin{equation}
\zxz{-1}{0}{0}{1} \matrixx{x}{y}{yD}{x}\zxz{-1}{0}{0}{1}=\matrixx{x}{-y}{-yD}{x}
\end{equation}
and thus
\begin{equation}
I(\varphi,\chi)=I\lrb{\pi\lrb{\zxz{-1}{0}{0}{1}}\varphi, \overline{\chi}}.
\end{equation}
So we shall always assume that $c(\theta{\overline{\chi}})\leq c(\theta\chi)$.

The first situation is when $c(\theta{\overline{\chi}})=0$, then the Tunnell-Saito's test requires $\theta{\overline{\chi}}$ to be trivial for $\WaldsI$ to be ever nonzero. In that case we can take $u=0$ and then there is a unique $v\mod \varpi^{\lceil n/2\rceil}$ such that $\WaldsI\neq 0$, and 
for this $v$ we have \begin{equation}
 \WaldsI\New{=\vol(\BF^\times\backslash\BE^\times)}=2.\yhb{better\ to\ add =\vol(\BF^\times\backslash\BE^\times)?}
\end{equation}

The second situation is when $0<c(\theta{\overline{\chi}})=2l\leq 2n$. In this case, $\alpha_{\theta{\overline{\chi}}}$ can be associated to $\theta\overline{\chi}$ by Lemma \ref{Lem:DualLiealgForChar} and $\varphi$ would be a test vector if $v$, $u$ are solutions of the following quadratic equation:
\begin{equation}\label{eq:new-necessary-same-ramified}
  \frac{D}{D'}v^2-\lrb{2\varpi^n\alpha_{\theta{\overline{\chi}}}\sqrt{D}-2\sqrt{\frac{D}{D'}}}v+(1-Du^2)\equiv 0\mod\varpi^{n-\lfloor \frac{l}{2}\rfloor}.
\end{equation}
This implies that for fixed $u$, the discriminant of the equation 
\begin{equation}
\Delta(u)=4\varpi^{n}\alpha_{\theta{\overline{\chi}}}\sqrt{D}\lrb{\varpi^{n}\alpha_{\theta{\overline{\chi}}}\sqrt{D}-2\sqrt{\frac{D}{D'}}}+4\frac{D}{D'}Du^2
\end{equation}
has to be a square $\mod \varpi^{n-\lfloor \frac{l}{2}\rfloor}$.
When $n-l$ is even, we can pick $u=0$ directly. Whether $\Delta(u)$ is a square is consistent with Tunnell-Saito's test.
When ${\Delta}(0)$ is indeed a square, we get two solutions of $v\mod \varpi^{ \lceil n/2\rceil}$. 
For each of these two solutions we have
\begin{equation}
\WaldsI=\frac{1}{q^{\lfloor l/2\rfloor}}.\yhb{maybe\ better\ refer\ to\ the\ Appendix\ for\ the\ result}
\end{equation}

Now if $n-l$ is odd, $v_\F({\Delta}(0))=n-l$ is odd, thus ${\Delta}(0)$ can never be a square. We need to pick $u$ such that $v_\F(u)=\frac{n-l-1}{2}$ and ${\Delta}(0)+4\frac{D}{D'}Du^2$ can be of higher evaluation and a square. Whether this is possible is again consistent with Tunnell-Saito's test.
 In this case it's possible to get more solutions of $v\mod\varpi^{\lceil n/2\rceil}$. For each solution we have
\begin{equation}
\WaldsI=\frac{1}{q^{\lfloor l/2\rfloor}}.\yhb{maybe\ better\ refer\ to\ the\ Appendix\ for\ the\ result}
\end{equation}

\section{Waldspurger's period integral using newforms}\label{newforms}
In the last section we reviewed the local Waldspurger's period integral for minimal vectors. In this section we show how to work out the local Waldspurger's period integral for newforms. The first step is to work out minimal vectors in the Kirillov model, which allows an explicit relation between minimal vectors and newforms. Then using bilinearity of Waldspurger's period integral, we can write the integral for newforms as a sum of integrals for minimal vectors. Using results from the last section, we show in Section \ref{Localcomputation} that there is a single diagonal term which is non-vanishing  for the special cases concerned in this paper. We further illustrate how to evaluate off-diagonal terms in Section \ref{Sec:testingInGeneral} in more general cases for possible future applications.

We need to deduce the local information from global characters in Section \ref{Localcomputation}, so we keep the subscripts to indicate local places there. In the other parts of this section we still omit subscript $v$.
\subsection{Kirillov model for minimal vectors}\label{Kirillov}
We first describe the minimal vectors explicitly in the Kirillov model. For this purpose, we choose a special shape for $B^1=U_\BL(1)K_{\fA_2}(2n+1)$ in the case $e_\BL=1$ and $c(\pi)=4n+2$. Recall we choose $D'$ such that  
\begin{equation}\label{Eq:specialAlphaTheta}
 \alpha_\theta= \frac{1}{\varpi_\BL^{c(\theta)}}\frac{1}{\sqrt{D'}} \mapsto \frac{1}{\varpi^{c(\theta)/e_\BL}}\matrixx{0}{\frac{1}{D'}}{1}{0}.
\end{equation}
We define the intertwining operator from $\pi$ to its Whittaker model by
\begin{equation}\label{eq:3.4:IntertwiningtoWhittaker}
\varphi \mapsto W_\varphi(g)=\int\limits_{\BF}\Phi_\varphi\lrb{\matrixx{\varpi^{\lfloor c(\pi)/2\rfloor}}{0}{0}{1}\matrixx{1}{n}{0}{1}g}\psi(-n)dn.
\end{equation}

\begin{lem}\label{lem:toricnewforminKirillov}
Up to a constant multiple, we have \New{the followings in the Kirillov model.} \yhb{add `in the Kirillov model'?}
\begin{enumerate}
\item When $c(\pi)=4n$, $\varphi_0=\Char(\varpi^{-2n}U_\BF(n)) $.
\item When $c(\pi)=2n+1$, $\varphi_0=\Char(\varpi^{-n}U_\BF(\lceil n/2 \rceil))$.
\item When $c(\pi)=4n+2$, $\varphi_0=\Char(\varpi^{-2n-1}U_\BF(n+1))$.
\end{enumerate}
\end{lem}
\begin{proof}
We only prove part (3) here. The first two results are similar and much easier. By the intertwining operator defined above, we have that
\begin{equation}
W_{\varphi_0}\lrb{\matrixx{a}{0}{0}{1}}=\int_\BF \Phi_{\varphi_0}\lrb{\matrixx{\varpi^{2n+1}a}{\varpi^{2n+1}n}{0}{1}}\psi(-n)dn.
\end{equation}
By the support of $\Phi_{\varphi_0}$ in Corollary \ref{Cor:MCofGeneralMinimalVec} we get that the integral is non-vanishing only when $a\in \varpi^{-2n-1}U_\BF(n+1)$ and $n\in \varpi^{-n-1}\CO_\BF$, in which case by \eqref{Eq:specialAlphaTheta} and definition of $\tilde{\theta}$,
\begin{equation}
W_{\varphi_0}\lrb{\matrixx{a}{0}{0}{1}}=\int_{n\in \varpi^{-n-1}\CO_\BF} \psi(\varpi^{-2n-1}\varpi^{2n+1}n)\psi(-n)dn=q^{n+1}.
\end{equation}
 So up to a constant, $\varphi_0$ in the Kirillov model is \New{$$\varphi_0(x)=W_{\varphi_0}\lrb{\matrixx{x}{0}{0}{1}}=\Char(\varpi^{-2n-1}U_\BF(n+1))(x),$$} \yhb{what is the relation between kirillov model and Whittaker model?}and the intertwining operator is not trivial. 
\end{proof}
From now on we no longer need the special shape of $\alpha_\theta$ or $B^1$. We do want to keep all test vectors $L^2-$normalised.
\begin{coro}\label{Cor:RelationNewMinimal}
The new form $\varphi_{new}$ can be related to $\varphi_0$ by the following formula
\begin{equation}
\varphi_{new}=\frac{1}{\sqrt{(q-1)q^{\lceil \frac{c(\theta)}{2 e_\BL}\rceil-1}}}\sum\limits_{x\in (\CO_\BF/\varpi^{\lceil \frac{c(\theta)}{2 e_\BL}\rceil} \CO_\BF)^\times} \pi\lrb{\matrixx{\varpi^{-c(\theta)/e_\BL}x}{0}{0}{1}}\varphi_0.
\end{equation}
Here $\varphi_0$ and $\varphi_{new}$ are both $L^2$-normalised.
\end{coro}
\begin{proof}
By the previous lemma one can uniformly write
\begin{equation}
\varphi_0=\sqrt{(q-1)q^{\lceil \frac{c(\theta)}{2 e_\BL}\rceil-1}}\Char(\varpi^{-c(\theta)/e_\BL} U_\BF(\lceil \frac{c(\theta)}{2 e_\BL}\rceil)).
\end{equation}
The coefficient comes from the $L^2$-normalisation of $\varphi_0$, as $$\Vol(U_\BF(\lceil \frac{c(\theta)}{2 e_\BL}\rceil))=\frac{1}{(q-1)q^{\lceil \frac{c(\theta)}{2 e_\BL}\rceil-1}}.$$
Then one just has to use that $\varphi_{new}=\Char(\CO_\BF^\times)$ in the Kirillov model. 
\end{proof}

To distinguish from $\beta(\varphi_1,\varphi_2)$ in \eqref{Eq:beta} which uses a different embedding for $K$, we denote
\begin{equation}\label{Eq:DefBetaI}
\BetaI{\varphi_1,\varphi_2}=\int\limits_{t\in \BF^\times\backslash\BE^\times}(\pi(t)\varphi_1,\varphi_2)\chi(t)dt
\end{equation}
for the embedding of $\BE$ as in \eqref{Eq:standardembedding}. Using its bilinearity and the previous result, we immediately have the following:

\begin{coro}\label{Cor:PeriodIntRelation}
Let $\widetilde{\varphi_{new}}=\pi\lrb{\matrixx{\varpi^{c(\theta)/e_\BL}}{0}{0}{1}}\varphi_{new}$, $\varphi_x=\pi\lrb{\matrixx{x}{0}{0}{1}}\varphi_0$ Then
\begin{equation}
\BetaI{\widetilde{\varphi_{new}},\widetilde{\varphi_{new}}}=\frac{1}{(q-1)q^{\lceil \frac{c(\theta)}{2 e_\BL}\rceil-1}}\sum\limits_{x,x'\in (\CO_\BF/\varpi^{\lceil \frac{c(\theta)}{2e_\BL}\rceil}\CO_\BF)^\times}\BetaI{\varphi_x,\varphi_{x'}}.
\end{equation}
\end{coro}
\New{
\begin{lem}\label{lem:Gl2-newform-crossterms}
\begin{enumerate}
\item Suppose that $\BetaI{\varphi_x,\varphi_x}=0$, then  $\BetaI{\varphi_{x'},\varphi_x}=\BetaI{\varphi_x,\varphi_{x'}}=0$ for any $x'$.\yhb{why? not trivial}
\item Suppose that $|\BetaI{\varphi_x,\varphi_x}|=|\BetaI{\varphi_{x'},\varphi_{x'}}|$. Then $|\BetaI{\varphi_{x},\varphi_{x'}}|= |\BetaI{\varphi_{x},\varphi_{x}}|$.
\end{enumerate}

\end{lem}
\begin{proof}
For any nontrivial functional $\mathcal{F}\in \Hom_{\BE^\times}(\pi\otimes \chi,\BC)$, we have $\BetaI{\varphi_1,\varphi_2}=C \mathcal{F}(\varphi_1)\overline{\mathcal{F}(\varphi_2)}$ for some non-zero constant $C$ independent of test vectors, as $\dim \Hom_{\BE^\times}(\pi\otimes \chi,\BC)\leq 1$.
Then
$$|\BetaI{\varphi_x,\varphi_x}|=|C\mathcal{F}(\varphi_x)^2|,$$
$$|\BetaI{\varphi_{x'},\varphi_{x'}}|=|C\mathcal{F}(\varphi_{x'})^2|,$$
$$|\BetaI{\varphi_{x},\varphi_{x'}}|=|C\mathcal{F}(\varphi_{x})\overline{\mathcal{F}(\varphi_{x'})}|.$$
Now the results are clear.
\end{proof}
}\yh{I moved this lemma from below and combined with zeroness}

The diagonal terms where $x=x'$ are already known by the last section. In general we need to understand the off-diagonal terms to evaluate the integral on newforms, but there are some special cases where only one diagonal term shows up and there are no off-diagonal terms as in Section \ref{Localcomputation}. 


\subsection{The special cases}\label{Localcomputation}
Now we specialise to the cases required in the proof of Theorem \ref{thm:GZ}. Recall that we are interested in the following Waldspurger's period integral
\begin{equation}
\beta^0_3(f_3, f_3)
=\int\limits_{t\in \BQ_3^\times\backslash K_3^\times} \Phi_{f_3}(t)\chi_3(t)dt
\end{equation}
 where $\Phi_{f_3}$ is the normalized matrix coefficient associated to $f_3$ and $K$ is embedded into $\GL_2(\BQ)$ as in Section \ref{heegner}. 
In particular we take $\varpi=3=q$, $D=-3$, $K\simeq\BE\simeq\BL\simeq\BQ(\sqrt{-3})$, $c(\theta_3)=c(\chi_3)=4$
. By Lemma \ref{lem:toricnewforminKirillov} we have for the minimal vector $\varphi_0=\Char(\varpi^{-2}U_\BF(1))$ in the Kirillov model. We shall give $\theta_3$ and $\chi_3$ explicitly using arithmetic information.
Note here that $K$ is embedded differently from the fixed embedding we have been using in \eqref{Eq:standardembedding}, but one can conveniently work out the relation between the two embeddings.

\subsubsection{Arithmetic information}
First of all we make use of the arithmetic information to give the local information explicitly. Recall that $K=\BQ(\sqrt{-3})$ is an imaginary quadratic field and $\CO_K=\BZ[\omega]$ is its ring of integers with $\omega=\frac{-1+\sqrt{-3}}{2}$. Let $\Theta:K^\times\backslash \BA_K^\times\ra \BC^\times$ be the unitary Hecke character associated to the CM elliptic curve ${E_9}_{/K}$.
Then $\Theta$ has conductor $9\CO_K$. For any place $v$ of $K$, let $\Theta_v$ be the local component of $\Theta$ at the place $v$. Then $\Theta_v$ is the character used to construct the local Weil-Deligne representation in Theorem \ref{LLC}. We denote $\Theta_3$ the 3-part of $\Theta$ and	denote $\pi_3$ be local representation of $\GL_2(\BQ_3)$ corresponding to $\Theta_3$. Note
\[\CO_{K,3}^\times/(1+9\CO_{K,3})\simeq \langle \pm 1\rangle ^{\BZ/2\BZ} \times\langle 1+\sqrt{-3}\rangle^{\BZ/3\BZ}\times\langle 1-\sqrt{-3}\rangle^{\BZ/3\BZ}\times \langle 1+3\sqrt{-3}\rangle^{\BZ/3\BZ}.\]
\begin{lem}\label{thetavalue}
The local character $\Theta_3$ is given explicitly by
\[\Theta_3(-1)=-1,\quad \Theta_3(1+\sqrt{-3})=\frac{-1-\sqrt{-3}}{2},\ \ \Theta_3(\sqrt{-3})=i,\]
\[\Theta_3(1-\sqrt{-3})=\frac{-1+\sqrt{-3}}{2},\quad \Theta_3(1+3\sqrt{-3})=\frac{-1+\sqrt{-3}}{2}.\]
\[\]
\end{lem}
\begin{proof}
It is well-known that $\Theta_\infty(x)=\frac{||x||}{x}$, (see for example \cite[Chapter II, Theorem 9.2]{Silvermanbook2} and note that we normalize it to make it unitary) and $\Theta$ is unramfied when $3\nmid v$. Note
\[\Theta_\infty(-1)\Theta_3(-1)=1,\quad \Theta_\infty(-1)=-1.\]
So $\Theta_3(-1)=-1$.

Let $\fp=(a)$ be a prime of $K$  relatively prime to $6$, with the unique generator $a\equiv 2\mod 3$.
By \cite[Chapter II, Example 10.6]{Silvermanbook2}, we have 
\[\Theta(\fp)=-\RN_\fp^{-1/2}\ov{\left(\frac{-3}{a}\right)}_6a.\]
Where $\left(\frac{\cdot}{a}\right)_6$ is the sixth power residue symbol and $\RN_\fp$ is the norm of $\fp$.
If $\fp=(5)$, then
\[\Theta_5(5)=-\ov{\left(\frac{-3}{5}\right)}_6=-1.\]
By
\[\Theta_\infty(10)\Theta_2(10)\Theta_3(10)\Theta_5(10)=1,\]
we get $\Theta_2(2)=-1.$ Since $\Theta_2$ is unramified and $1+\sqrt{-3}$ is another uniformizer of $(2)$, we see $\Theta_2(1\pm\sqrt{-3})=\Theta_2(2)=-1$.
By
\[\Theta_\infty(1\pm\sqrt{-3})\Theta_2(1\pm\sqrt{-3})\Theta_3(1\pm\sqrt{-3})=1,\quad \Theta_\infty(1\pm\sqrt{-3})=\frac{2}{1\pm\sqrt{-3}},\]
we get
\[\Theta_3(1\pm\sqrt{-3})=\frac{-1\mp\sqrt{-3}}{2}.\]

Let $v,\ov{v}$ be the places above $7$ with corresponding prime ideals $\fp_v=\left(\frac{1+3\sqrt{-3}}{2}\right)$ and $\fp_{\ov{v}}=\left(\frac{1-3\sqrt{-3}}{2}\right)$. Then
\[\Theta_v(1+3\sqrt{-3})=-\sqrt{7}^{-1}\cdot \ov{\left(\frac{-3}{\frac{1+3\sqrt{-3}}{2}}\right)}_6\cdot \frac{1+3\sqrt{-3}}{2}=-\sqrt{7}^{-1}\cdot\frac{-1-\sqrt{-3}}{2}\cdot \frac{1+3\sqrt{-3}}{2}.\]
By
\[\Theta_\infty(1+3\sqrt{-3})\Theta_2(1+3\sqrt{-3})\Theta_3(1+3\sqrt{-3})\Theta_v(1+3\sqrt{-3})=1,\quad \Theta_\infty(1+3\sqrt{-3})=\frac{2\sqrt{7}}{1+3\sqrt{-3}},\]
we get
\[\Theta_3(1+3\sqrt{-3})=\frac{-1+\sqrt{-3}}{2}.\]

At last, since $\Theta_\infty\Theta_3(\sqrt{-3})=1$ and $\Theta_{\infty}(\sqrt{-3})=-i$,
we have $\Theta_3(\sqrt{-3})=i$.
\end{proof}

Let $p\equiv 4,7\mod 9$ be a prime.  Let $\chi:G_K\ra \CO_K^\times$ be the character given by $\chi(\sigma)=(\sqrt[3]{3p})^{\sigma-1}$. We also view $\chi$ as a Hecke character on $\BA_K$ by the Artin map.
The local  character  $\chi_3$ has conductor $\BZ_3^\times(1+9\CO_{K,3})$, and hence it is in fact a character of the quotient group $\CO_{K,3}^\times/\BZ_3^\times(1+9\CO_{K,3})$. Note that
\[\CO_{K,3}^\times/\BZ_3^\times(1+9\CO_{K,3})\simeq \langle 1+\sqrt{-3}\rangle^{\BZ/3\BZ}\times\langle 1+3\sqrt{-3}\rangle^{\BZ/3\BZ}.\]
We have the following
\begin{lem} \label{chi}
The local character $\chi_3$ is given explicitly by the following table:
\begin{center}
\begin{tabular}{|c|c|c|c|c|c|}
\hline
$p\mod 9$&$\chi_3(1+\sqrt{-3})$&$\chi_3(1+3\sqrt{-3})$&$\chi_3(\sqrt{-3})$\\
\hline
$4$&$\omega$&$\omega$&$1$\\
\hline
$7$&$\omega^2$&$\omega$&$1$\\
\hline
\end{tabular}
\end{center}
\end{lem}

\begin{proof}
This follows directly from the explicit local class field theory, see the proof of Proposition \ref{LCF}.
We just remark that $1+9\CO_{K,3}\subset K_3^{\times 3}$, which follows easily by considering the $3$-adic valuation of the terms of the binomial expansion:
\begin{equation}\label{newton}
(1+9x)^{\frac{1}{3}}=1+\sum_{n\geq 1} \frac{1/3(1/3-1)\cdots(1/3-n+1)}{n!}(9x)^n
\end{equation}
where $x\in \CO_{K,3}$. The $3$-adic valuation of the $n$-th term
\[\frac{1/3(1/3-1)\cdots(1/3-n+1)}{n!}(9x)^n\]
is greater than $n/3$ for sufficiently large $n$. Then the right-hand side of (\ref{newton}) is convergent, and we see $1+9\CO_{K,3}\subset K_3^{\times 3}$.

Hence for any $t\in K_3^\times$,
\[\chi_3(t)=\left(\sqrt[3]{3p}\right)^{\sigma_t-1}=\lrb{\frac{t,3p}{K_3;3}}=\left\{\begin{aligned}\lrb{\frac{t,12}{K_3;3}},&\quad p\equiv 4\mod 9;\\ \lrb{\frac{t,21}{K_3;3}},&\quad p\equiv 7\mod 9.\\\end{aligned}\right.\]
Recall $\sigma_t$ is the  image of $t$ under the the Artin map, and $\lrb{\frac{\cdot,\cdot}{K_3;3}}$ denotes the $3$-rd Hilbert symbol over $K_3^\times$ as before. Using the local and global principal, it is straight-forward to compute the values of $\chi_3$ as in the above table.
\end{proof}

\subsubsection{Local period integral}

From Section \ref{Sec:WaldsForMinimalVec}, we see that the test vector issue for Waldspurger's period integral is closely related to $c(\theta_3\ov\chi_3)$ or $c(\theta_3\chi_3)$ and some further details like $\alpha_{\theta_3\ov\chi_3}$. We can work out these  by using Lemma \ref{thetavalue}, \ref{chi}, and the relation between $\theta_3$ and $\Theta_3$ as in Theorem \ref{LLC}.
\begin{coro}
\label{thetachivalue}
If $p\equiv 4\mod 9$,
the local character $\Theta_3\ov\chi_3$ is given explicitly by
\[\Theta_3\ov\chi_3(-1)=-1,\quad \Theta_3\ov\chi_3(1+\sqrt{-3})=\omega,\]
\[\Theta_3\ov\chi_3(1-\sqrt{-3})=\omega^2,\quad \Theta_3\ov\chi_3(1+3\sqrt{-3})=1,\quad \Theta_3\ov\chi_3(\sqrt{-3})=i.\]

If $p\equiv 7\mod 9$,
the local character $\Theta_3\ov\chi_3$ is given explicitly by
\[\Theta_3\ov\chi_3(-1)=-1,\quad \Theta_3\ov\chi_3(1+\sqrt{-3})=1,\]
\[\Theta_3\ov\chi_3(1-\sqrt{-3})=1,\quad \Theta_3\ov\chi_3(1+3\sqrt{-3})=1,\quad \Theta_3\ov\chi_3(\sqrt{-3})=i.\]
\end{coro}
Now we can prove the following key Lemma in our special case.
\begin{lem}\label{Cor:AllnecessaryFormulationForThetaChi}
When $p\equiv 7\mod 9$, we have $\theta_3\ov\chi_3=1$. When $p\equiv 4\mod 9$, we have $c(\theta_3\ov\chi_3)=2$ and $\alpha_{\theta_3\ov\chi_3}=\frac{1}{3\sqrt{-3}}$.
\end{lem}
\begin{proof}
Let $\psi_3$ be the additive character such that $\psi_3(x)=e^{2\pi i \iota(x)}$ where $\iota:\BQ_3\rightarrow \BQ_3/\BZ_3 \subset \BQ/\BZ$ is the map given by $x\mapsto -x\mod \BZ_3$ which is compatible with the choice in \cite{CST14}. Let $\psi_{K_3}(x)=\psi_3\circ \Tr_{K_3/\BQ_3}(x)$, be the additive character of $K_3$. 

$\alpha_{\Theta_3}$ is the number associated to $\Theta_3$ as in Lemma \ref{Lem:DualLiealgForChar} so that
\[\Theta_3(1+x)=\psi_{K_3}(\alpha_{\Theta_3} x), \]
for any $x$ satisfying $v_{K_3}(x)\geq c(\Theta_3)/2=2$. By the definition of $\psi_{K_3}$ and Lemma \ref{thetavalue}, we know that $\alpha_{\Theta_3}=\frac{1}{9\sqrt{-3}}$.  Now let $\eta_3$ be the quadratic character associated to the quadratic field extension $K_3/\BQ_3$.
Then by \cite[Proposition 34.3]{BushnellHenniart:06a}, $\lambda_{K_3/\BQ_3}(\psi')=\tau(\eta_3,\psi_3')/\sqrt{3}=-i$, here $\tau(\eta_3,\psi_3')$ is the Gauss sum and $\psi_3'(x)=\psi_3(\frac{x}{3})$ is the additive character of level one. By \cite[Lemma 5.1]{Langlands}, $\lambda_{K_3/\BQ_3}(\psi_3)=\eta_3(3)\lambda_{K_3/\BQ_3}(\psi_3')=-i$. Then $\Delta_{\theta_3}$ is the unique level one character of $K_3$ such that $\Delta_{\theta_3}|_{\BZ_3^\times}=\eta_3$ and
\[\Delta_{\theta_3}(\sqrt{-3})=\eta((\sqrt{-3})^3\alpha_{\Theta_3})\lambda_{K_3/\BQ_3}(\psi_3)^{3}=-i.\]
Recall that $\theta_3=\Theta_3\Delta_{\theta_3}$. Then by Corollary \ref{thetachivalue} we can easily check that: 

\begin{enumerate}
\item If $p\equiv 7\mod 9$, $\theta_3\ov\chi_3$ is the trivial character. 
\item If $p\equiv 4\mod 9$, $\theta_3\ov\chi_3$ is of level 2 and by definition we can choose $\alpha_{\theta_3\ov\chi_3}=\frac{1}{3\sqrt{-3}}$.
\end{enumerate}
\end{proof}

Recall that we embed $K$ into $\rm{M}_2(\BQ)$ by linearly extending the following map in section \ref{heegner}. 
\begin{equation}\label{matrix}
\sqrt{-3}\mapsto \matrixx{4p+17+72/p}{-8p/9-4-18/p}{18p+72+288/p}{-4p-17-72/p}=:\matrixx{a}{3^{-2}b}{3^3c}{-a}
\end{equation}
with $3||a$ if $p\equiv 4\mod 9$, $9||a$ if $p\equiv 7\mod 9$, $b\equiv p\mod 9$ and $c\equiv -1\mod 9$.  Then $\Nm(\sqrt{-3})=-a^2-3bc=3$.
 
\begin{prop}\label{Prop:TestingForNew}
Suppose $\Vol(\BZ_3^\times\backslash\CO_{K,3}^\times)=1$ so that $\Vol(\BQ_3^\times\backslash K_3^\times)=2$.
For $f_3$ being the standard $L^2-$normalised newform, $K$ being embedded as in \eqref{matrix} and $\theta_3,\chi_3$ as given above, we have
\begin{equation}
 \beta_3(f_3,f_3)=\begin{cases}
 1, &\text{\ if }p\equiv 7\mod 9\\
 1/2, &\text{\ if }p\equiv 4\mod 9.
 \end{cases}
\end{equation}
\end{prop}

\begin{proof}
To evaluate $f_3$ for the embedding in \eqref{matrix} is equivalent to use the standard embedding \eqref{Eq:standardembedding} of $\BE$ and use different translate of the newform.
In particular the embedding in \eqref{matrix} is conjugate to the standard embedding by the following. 
\begin{equation}
    \matrixx{a}{3^{-2}b}{3^3c}{-a}=\matrixx{-9c}{a/3}{0}{1}^{-1}\matrixx{0}{1}{D}{0}\matrixx{-9c}{a/3}{0}{1}.
\end{equation}
Thus  we have
\begin{align}
\beta_3(f_3,f_3)&=\int\limits_{\BF^\times\backslash \BE^\times}(\pi_3\lrb{\matrixx{-9c}{a/3}{0}{1}^{-1}t\matrixx{-9c}{a/3}{0}{1}}f_3,f_3)\chi(t)dt\\
&=\int\limits_{\BF^\times\backslash \BE^\times}(\pi_3\lrb{t\matrixx{-9c}{a/3}{0}{1}}f_3,\pi_3\lrb{\matrixx{-9c}{a/3}{0}{1}f_3})\chi(t)dt\notag,
\end{align}
which is by definition $\BetaI{\pi_3\lrb{\matrixx{-9c}{a/3}{0}{1}}f_3,\pi_3\lrb{\matrixx{-9c}{a/3}{0}{1}f_3}}$ for the bilinear pairing as in \eqref{Eq:DefBetaI} and the standard embedding as in \eqref{Eq:standardembedding}.
Note that by Corollary \ref{Cor:RelationNewMinimal}
$$\pi_3\lrb{\matrixx{-9c}{a/3}{0}{1}}f_3
=\frac{1}{\sqrt{2}}\sum\limits_{x\in (O_v/\varpi O_v)^\times}\pi_3\lrb{\matrixx{1}{a/3}{0}{1}\matrixx{x}{0}{0}{1}}\varphi_0$$
where $\varphi_0$ is the minimal test vector.

Now there are two cases. When $p\equiv 7\mod 9$, $9|| a$ and the action of $\matrixx{1}{a/3}{0}{1}$ on $\varphi_x=\pi_3\lrb{\matrixx{x}{0}{0}{1}}\varphi_0$ is by a simple character. By the $l=0$ case in Section \ref{Sec:WaldsForMinimalVec}, we have a unique $x\mod\varpi$ for which $\BetaI{\varphi_x,\varphi_x}$ is nonvanishing. So there are no off-diagonal terms in Corollary \ref{Cor:PeriodIntRelation}, and we have
\begin{equation}
\beta_3\lrb{f_3,f_3 }=\frac{1}{(q-1)q^{\lceil \frac{c(\theta)}{2 e_\BL}\rceil-1}} \BetaI{\varphi_x,\varphi_x}=\frac{1}{2}\cdot 2=1.
\end{equation}

When $p\equiv 4\mod 9$, we have $3||a$ and $u=a/3$. By Lemma \ref{Cor:AllnecessaryFormulationForThetaChi}, this is the case $l=1$ and $n-l=1$ is odd. By the choice in Section \ref{Sec:WaldsForMinimalVec}, 
$$D'=\frac{1}{\alpha_\theta^2\varpi_\BL^{2c(\theta)}}=-3.$$

By Lemma \ref{Cor:AllnecessaryFormulationForThetaChi}, $\alpha_{\theta_3\chi_3^{-1}}=\frac{1}{3\sqrt{-3}}$ in this case, and we have 
\begin{align}
\Delta(u)&= 4\varpi^{n}\alpha_{\theta_3\ov\chi_3}\sqrt{D}\lrb{\varpi^{n}\alpha_{\theta_3\ov\chi_3}\sqrt{D}-2\sqrt{\frac{D}{D'}}}+4\frac{D}{D'}Du^2 \\
&\equiv 4\cdot 9\cdot \frac{1}{3\sqrt{-3}} \cdot \sqrt{-3}\cdot (-2)+4 \cdot (-3)\frac{a^2}{9}\mod{\varpi^2}\notag\\
&\equiv  -8\cdot 3-4\cdot 3 \mod{\varpi^2}\notag\\
&\equiv 0 \mod{\varpi^2}\notag.
\end{align}
$\Delta(u)$ is indeed  congruent to a square. So we also get a unique solution of $x\mod \varpi$ from \eqref{eq:new-necessary-same-ramified}, and 
\begin{equation}
 \beta_3\lrb{f_3,f_3}=\frac{1}{(q-1)q^{\lceil \frac{c(\theta)}{2 e_\BL}\rceil-1}}\frac{1}{q^{\lfloor l/2\rfloor}}=\frac{1}{2}.
\end{equation}
\end{proof}

\begin{coro}\label{ration}
For the admissible test vector $f'_3$ and the standard new form $f_3$ we have 
\[\frac{\beta_3^0(f'_3,f'_3)}{\beta_3^0(f_3,f_3)}=\begin{cases}
 2, &\text{\ if }p\equiv 7\mod 9,\\
 4, &\text{\ if }p\equiv 4\mod 9.
 \end{cases}\]
\end{coro}
\begin{proof}
Keep the normalization of the volumes in Proposition \ref{Prop:TestingForNew}. By definition of $f'$, we have $\beta_3^0(f'_3,f'_3)=\Vol(\BQ_3^\times\backslash K_3^\times)=2$. Then the corollary follows from Proposition \ref{Prop:TestingForNew}, considering that $f_3$ is $L^2$-normalised there.
\end{proof}

\subsection{The general cases}\label{Sec:testingInGeneral}
We see from the discussion above that when there is a unique solution of $v\mod \varpi^{\lceil n/2\rceil}$  for \eqref{eq:new-necessary-same-ramified} with fixed $u$, one can easily predict the local Waldspurger's period integral for newforms from the corresponding integral for minimal vectors. Here we explain how to work out the off-diagonal terms in Corollary \ref{Cor:PeriodIntRelation} in more general cases. It is not necessary for the main purpose of this paper, but may be useful for future applications.

By the support of the local Waldspurger's period integral $\WaldsI$, we mean the subset $\BE^\times\cap \Supp \Phi_{\varphi}$.
The main idea is to get the size of the off-diagonal terms in Lemma \ref{lem:Gl2-newform-crossterms} by purely representation theoretical consideration, and to get the support of the integral in Lemma \ref{Lem:offdiagSupp}. The volume of the support of the integral is exactly the size of the integral, while the integrand is absolutely bounded by 1. This forces the integrand to be constant (with absolute value 1) on the support of the integral. Then one can easily detect this constant by looking at the value of the integrand at any point in the support of the integral.

For simplicity, however, we stay in the setting where $\BE\simeq \BL$ are ramified, $0<c(\theta{\overline{\chi}})=2l\leq 2n$. We further assume that $n-l$ is even. By Section \ref{Sec:WaldsForMinimalVec}, we can pick $u=0$, and there exists 2 solutions $v,v' \mod \varpi^{\lceil n/2\rceil}$ to \eqref{eq:new-necessary-same-ramified}, while the local period integral is always $\frac{1}{q^{\lfloor l/2\rfloor}}$ for each solution.

\yh{I've moved a lemma away from here}
According to Lemma \ref{lem:Gl2-newform-crossterms}, we can write that $\BetaI{\varphi_{v},\varphi_{v'}}= \gamma \BetaI{\varphi_{v},\varphi_{v}}$ for some $\gamma$ with $|\gamma|=1$. Then
\begin{align}
I(\varphi_{new},\chi)&=\frac{1}{(q-1)q^{\lceil \frac{c(\theta)}{2 e_\BL}\rceil-1}}(\BetaI{\varphi_{v},\varphi_{v}}+\BetaI{\varphi_{v},\varphi_{v'}}+\BetaI{\varphi_{v'},\varphi_{v}}+\BetaI{\varphi_{v'},\varphi_{v'}})\\
&=\frac{1}{(q-1)q^{\lceil \frac{c(\theta)}{2 e_\BL}\rceil-1}}\frac{1}{q^{\lfloor l/2\rfloor}}(1+\gamma)(1+\overline{\gamma}). \notag
\end{align}

To study $\gamma$,  we first study the support of the integral. Without loss of generality, we can assume that $v$, $v'$ satisfy
\begin{equation}\label{Eq:exactEqvv'}
 \frac{D}{D'}v^2-\lrb{2\varpi^n\alpha_{\theta{\overline{\chi}}}\sqrt{D}-2\sqrt{\frac{D}{D'}}}v+1=0,
\end{equation}
compared to \eqref{eq:new-necessary-same-ramified} while taking $u=0$.

Let $k=\matrixx{v}{0}{0}{1}$ and $k'=\matrixx{v'}{0}{0}{1}$.
Then for $t=\matrixx{a}{b}{bD}{a}$,
$$k'^{-1}tk=\matrixx{v v'^{-1}a }{v'^{-1}b}{vbD}{a},$$
and
\begin{equation}
\BetaI{\varphi_{v},\varphi_{v'}}=\int\limits_{\BF^\times\backslash\BE^\times} \Phi_{\varphi_0}(k'^{-1}tk)\chi(t)dt.
\end{equation}
\begin{lem}\label{Lem:offdiagSupp}
For the setting as above, we have $vv'D=D'$, and $v_\F(\frac{v}{v'}-1)= \frac{n-l}{2}$.
In particular the support of the integral 
is $v_\F(b)=0, v_\F(a)\geq  \lceil\frac{l+1}{2}\rceil$.
\end{lem}
\begin{proof}
According to \eqref{Eq:exactEqvv'}, $v$ and $v'$ satisfy 
$$v v'=\frac{D'}{D},\text{ \ \ } v+v'=2(\varpi^{n}\alpha_{\theta{\overline{\chi}}}\sqrt{D'}-1)\sqrt{\frac{D'}{D}}$$
so the first result is direct.
For the second result, we note that
\begin{equation}
\lrb{\frac{v}{v'}-1}^2=\frac{(v+v')^2-4vv'}{v'^2}=\frac{D'}{D}\frac{4\varpi^{n}\alpha_{\theta{\overline{\chi}}}\sqrt{D'}(\varpi^{n}\alpha_{\theta{\overline{\chi}}}\sqrt{D'}-2)}{v'^2}.
\end{equation}
Thus $v_\F((\frac{v}{v'}-1)^2)=n-l$ and $v_\F(\frac{v}{v'}-1)= \frac{n-l}{2}$. 
Now for $k'^{-1}tk\in J=\BL^\times K_{\fA_{e_\BL}}(n)$, there are two parts to consider: either  $v_\F(b)=0$, $v_\F(a)>0$ or $v_\F(a)=0$, $v_\F(b)\geq 0$. 
In the first case, since $$\varpi_\BL \CB_{e_\BL} ^n=\matrixx{\varpi^{ \lceil \frac{n+1}{2}\rceil}\CO_\BF}{\varpi^{ \lfloor \frac{n+1}{2}\rfloor}\CO_\BF}{\varpi^{ \lfloor \frac{n+1}{2}\rfloor+1}\CO_\BF}{\varpi^{ \lceil \frac{n+1}{2}\rceil}\CO_\BF}$$ we must have
\begin{equation}
a(\frac{v}{v'}-1)\equiv 0 \mod \varpi^{\lceil \frac{n+1}{2}\rceil}, \text{\ \ }vbD-v'^{-1}bD'\equiv 0 \mod \varpi^{\lfloor\frac{n+1}{2}\rfloor +1}.
\end{equation}
The second equation is automatic as $v v'=\frac{D'}{D}$. From the first equation and the computation for $v_\F(\frac{v}{v'}-1)$ above, we get $v_\F(a)\geq \lceil\frac{l+1}{2}\rceil$.
Using similar argument, one can easily show that it is impossible for $k'^{-1}tk\in \BL^\times K_{\fA_{e_\BL}}(n)$ when $v_\F(a)=0$, $v_\F(b)\geq 0$.
\end{proof}
\begin{prop}\label{Prop:TestingonNewform}
$$I(\varphi_{new},\chi)=\frac{1}{(q-1)q^{\lceil \frac{c(\theta)}{2 e_\BL}\rceil-1}}\frac{1}{q^{\lfloor l/2\rfloor}} (1+\theta\chi(\sqrt{D}))^2    .$$
In particular it is either $0$ or asymptotically $\frac{1}{C(\pi\times\pi_{\chi})^{1/4}}$.
\end{prop}
\begin{proof}
We already know that $|\BetaI{\varphi_{v},\varphi_{v'}}|=\frac{1}{q^{\lfloor l/2\rfloor}} $. By the previous lemma the support of the integral also has volume $\frac{1}{q^{\lfloor l/2\rfloor}}$, while the integrand  $|<\pi(t)\varphi_v,\varphi_{v'}>\chi(t)|\leq 1$. Then  $<\pi(t)\varphi_v,\varphi_{v'}>\chi(t)$ must be some constant $\gamma$ on the whole support with $|\gamma|=1$. To detect this constant we just have to take $t=\matrixx{0}{1}{D}{0}$. Then
\begin{equation}
\gamma=\Phi_{\varphi_0}\lrb{\matrixx{0}{v'^{-1}}{vD}{0}}\chi(\sqrt{D})=\theta(v'^{-1}\sqrt{D'})\chi(\sqrt{D})=\theta\chi(\sqrt{D}).
\end{equation}
In the last equality we have used $\theta|_{\BF^\times}=1$ so that $$\theta(v'^{-1}\sqrt{D'})=\theta\lrb{v'^{-1}\frac{\sqrt{D'}}{\sqrt{D}}\sqrt{D}}=\theta(\sqrt{D}).$$ Note that $\theta\chi(\sqrt{D})=\pm 1$.
The last statement is easy to check.
\end{proof}
\begin{remark}
The strategy should work in full generality.
\end{remark}

\appendix
\section{Explicit Waldspurger's period integral using minimal vectors}
As \cite{HN18} is still a preprint, we take part of it as appendix here which is closely related to this paper. If \cite{HN18} come out, we can omit this appendix conveniently.\yhb{I add this sentence} This appendix is purely local and we omit all subscripts $v$.

\subsection{Explicit Tunnell-Saito's $\epsilon-$value test}\label{Sec:epsilontest}
Here we recall the explicit results of Tunnell when $p\neq 2$, and later on we shall show that using minimal vectors for Waldspurger's period integral can effectively reprove Tunnell's work.
\begin{thm}[\cite{Tunnell:83a} \cite{Saito:93a}]\label{Tunnell}
Suppose that $w_\pi\cdot\chi|_{\F^\times}=1$.
The space $\Hom_{\E^\times}(\pi^\B\otimes {\chi},\C)$ is at most one-dimensional. It is nonzero if and only if 
\begin{equation}
\epsilon(\pi_\E\times{\chi})=\chi(-1)\epsilon(\B).
\end{equation}
Here $\pi_\E$ is the base change of $\pi$ to $\E$. $\epsilon(\B)=1$ if it is a matrix algebra, and $-1$ if it's a division algebra.
\end{thm}

Let $\pi$ be associated to $\sigma=\Ind_{\L}^{\F}\Theta$ for a character $\Theta $ over $\L$ with $c(\Theta)=n$. Let $\chi$ be a character over $\E$ with $c(\chi)=m$.
As noted before, we assume $p\neq 2$ and trivial central character. Further we restrict ourselves to the case $c(\pi)\geq c(\pi_{\chi})$ (otherwise it's always on $\GL_2$ side).
Tunnell's work actually give the $\epsilon-$value explicitly when $p\neq 2$. For the cases concerned in this paper, we take the following lemma directly from  \cite{Tunnell:83a}.

\begin{lem}\label{lem:Epsilontestsameramified}
Let $\E$ and $\L$ be the same ramified quadratic extensions with uniformizer $\varpi_\E$ such that $\varpi_\E^2=\varpi_\F$. 
Let $\eta$ be a character of $\E^\times$ extending the quadratic character $\eta_{\E/\F}$ of $\F^\times$. 
Then $\epsilon(\pi_\E\times{\chi})=-1$ iff one of the followings is true
\begin{enumerate}
\item $c(\Theta{\chi}\eta)=c(\overline{\Theta}{\chi}\eta)=c(\Theta)$ and, for all $x\in \CO_\F$, $\chi(1+\varpi_\E^jx)=\Theta(1+\varpi_\E^j \delta x)$ where $j=c(\Theta)-1$ and $\delta \in \CO_\F/\varpi \CO_\F$ satisfies $\delta ^2-1$ is not a square in $\CO_\F/\varpi \CO_\F$.
\item $0<c(\Theta{\chi}\eta)<c(\Theta)$ and the character $\nu=\Theta{\chi}\eta$ satisfies $\nu(1+\varpi_\E^{c(\nu)-1}x)=\Theta(1+\varpi_\E^{c(\Theta)-1} \delta x)$ for $x\in \CO_\F$ where $\delta \in \CO_\F^\times$ satisfies $-2\delta (-1)^{(c(\nu)+c(\Theta))/2}$ is not a square $\mod \varpi \CO_\F$.
\item $0<c(\overline{\Theta}{\chi}\eta)<c(\Theta)$ and the character $\nu=\overline{\Theta}{\chi}\eta$ satisfies $\nu(1+\varpi_\E^{c(\nu)-1}x)=\Theta(1+\varpi_\E^{c(\Theta)-1} \delta x)$ for $x\in \CO_\F$ where $\delta \in \CO_\F^\times$ satisfies $2\delta  (-1)^{(c(\nu)+c(\Theta))/2}$ is not a square $\mod \varpi \CO_\F$.
\item $c(\nu)=0$ for $\nu=\Theta{\chi}\eta$ or $\overline{\Theta}{\chi}\eta$, and $\nu(\varpi_\E)=-1$.
\end{enumerate}

\end{lem}
Note that we don't need to worry about the case when $c(\nu)=1$. This is because when $\nu |_{\F^\times}=1$ and $\E$ is ramified, such characters don't exist.

\subsection{Testing Waldspurger with minimal vector}

 For simplicity we shall focus on the case when $\BL=\BE$ are ramified,  $\pi$ has trivial central character, $c(\pi)=2n+1$ is odd and $c(\pi_{\chi})\leq c(\pi)$.
 Let $\varphi_0$ be a minimal vector with matrix coefficient  described  as in Corollary \ref{Cor:MCofGeneralMinimalVec}. Let $$\varphi=\pi\lrb{\zxz{1}{u}{0}{1}\zxz{v}{0}{0}{1}}\varphi_0$$
for some $u,v\in \F^\times$.
Denote $k=\zxz{1}{u}{0}{1}\zxz{v}{0}{0}{1}$. Then for $\varphi$ we have
\begin{align}
I(\varphi,\chi)&=\int\limits_{\F^\times\backslash\E^\times}\Phi_\varphi(t){\chi}(t)dt\\
&=\int\limits_{\F^\times\backslash\E^\times}\Phi_{\varphi_0}(k^{-1}tk){\chi}(t)dt \notag
\end{align}
The idea for this integral is very simple. When $\L$ is ramified, $\Phi_{\varphi_0}$ is a multiplicative character on the support. 
If there exists some $t_0\in \E^\times\cap \text{Supp}\Phi_{\varphi}$ such that $\Phi_{\varphi_0}(k^{-1}t_0k){\chi}(t_0)\neq  1$, one can make a change of variable and see that the integral must be vanishing.
Then the integral is nonvanishing iff $\Phi_{\varphi_0}(k^{-1}tk){\chi}(t)=1$ on $\E^\times\cap \text{Supp}\Phi_{\varphi}$. 
In that case the local integral is simply the volume of  $\E^\times\cap \text{Supp}\Phi_{\varphi}$.

From now on we fix the embeddings of $\L$ and $\E$ as in Section \ref{Sec:WaldsForMinimalVec}.

 We shall also assume that $$v_\F(v)=0\text{\  and }v_\F(u)\geq 0.$$ In principle we need to consider all possible valuations to cover all possible test vectors. But it turns out the test vectors with these restrictions already suffice.

Let $t=\zxz{a}{b}{bD}{a}\in \E$, then
\begin{equation}\label{Eq:offdiagonalembedding}
 k^{-1}tk=\zxz{a-bDu}{v^{-1}b(1-Du^2)}{vbD}{a+bDu}.
\end{equation}

Also recall that

$$\mathcal{B}^n=\zxz{\varpi^{\lceil n/2\rceil}\CO_\F}{\varpi^{\lfloor n/2\rfloor}\CO_\F}{\varpi^{\lfloor n/2\rfloor+1}\CO_\F}{\varpi^{\lceil n/2\rceil}\CO_\F}.$$

\subsubsection{Case  $0<c(\theta{\overline{\chi}})=2l\leq 2n=c(\theta\chi)$}
Since $\chi|_{\F^\times}$ and $\theta|_{\F^\times}$ are both trivial, $c(\theta{\overline{\chi}})$ must be even.  The case $c(\theta\chi)\leq 2n$ is parallel, so we shall skip this case here.

For $t=\zxz{a}{b}{bD}{a}$  with $v_\F(b)\geq v_\F(a)=0$, $k^{-1}tk\in \L^\times K_\mathfrak{A}(n)$ if and only if
\begin{equation}
a-bDu\equiv a+bDu\mod{\varpi^{\lceil n/2\rceil}},
\end{equation}
and
\begin{equation}
vb\frac{D}{D'}\equiv v^{-1}b(1-Du^2) \mod{\varpi^{\lfloor n/2\rfloor}}. 
\end{equation}
By simple manipulations, they are equivalent to that
\begin{equation}\label{Eq:Support1}
bDu\equiv 0\mod{\varpi^{\lceil n/2\rceil}},
\end{equation}
and
\begin{equation}\label{Eq:Supp2}
vb\frac{D}{D'}\equiv v^{-1}b \mod{\varpi^{\lfloor n/2\rfloor}}.
\end{equation}
Note that $\frac{D}{D'}$ is indeed a square. 
When $v_\F(a)>v_\F(b)=0$, note that $\varpi_\L \mathcal{B}_{e_\L}^n=\mathcal{B}_{e_\L}^{n+1}$. Then $k^{-1}tk\in \L^\times K_\mathfrak{A}(n)$ if and only if
\begin{equation}
bDu\equiv 0 \mod\varpi^{\lfloor n/2\rfloor+1},
\end{equation}
\begin{equation}
vb\frac{D}{D'}\equiv v^{-1}b\mod\varpi^{\lceil n/2\rceil}.
\end{equation}
They are true  only if  
\begin{equation}\label{eq:addtional-support1}
Du\equiv 0\mod\varpi^{\lfloor n/2\rfloor+1},
\end{equation}
and
\begin{equation}\label{eq:additional-support2}
v\frac{D}{D'}\equiv v^{-1}\mod\varpi^{\lceil n/2\rceil}
\end{equation}  
as $v_\F(b)=0$.

Note here that \eqref{eq:addtional-support1} and \eqref{eq:additional-support2} are respectively stronger congruence conditions than \eqref{Eq:Support1} and \eqref{Eq:Supp2}. This implies that if the support of the integral contains some of the part $v_\F(a)>v_\F(b)=0$, it actually contains the whole torus, and the local period integral is nonvanishing iff $2l=c(\theta\overline{\chi})=0$ (or $c(\theta\chi)=0$) which is to be discussed in the next subsection. By considering the case $l>0$ we can assume that the support of the integral contains only the part $v_\F(b)\geq v_\F(a)=0$.
Further more from \eqref{Eq:Supp2}, we get that
\begin{equation}
 v^{-1}b\equiv \pm\sqrt{\frac{D}{D'}}b\mod\varpi^{\lfloor n/2\rfloor}.
\end{equation}
The case $ v^{-1}b\equiv \sqrt{\frac{D}{D'}}b$ is parallel to the case $ v^{-1}b\equiv -\sqrt{\frac{D}{D'}}b$ and is related to the situation $c(\theta\chi)<2n$.

In the following lemma we explain how to write down the value of matrix coefficient, once we know $k^{-1}tk\in \L^\times K_\mathfrak{A}(n)$.
\begin{lem}\label{lem:supportimpliesclosetorus}
Suppose that $v_\F(b)\geq v_\F(a)=0$, $bDu\equiv 0\mod\varpi^{\lceil n/2\rceil}$  and $v^{-1}b\equiv -\sqrt{\frac{D}{D'}}b\mod\varpi^{\lfloor n/2\rfloor}$, then we have
\begin{align}\label{Eq:rewritingconjugation}
k^{-1}t k&=\frac{1}{a^2-b^2D}\zxz{a}{-b\sqrt{\frac{D}{D'}}}{-b\sqrt{DD'}}{a}\\
&\times\zxz{a^2-abDu+vb^2\sqrt{\frac{D}{D'}}D}
{ab(v^{-1}+\sqrt{\frac{D}{D'}})-bDu(v^{-1}au-b\sqrt{\frac{D}{D'}})}
{abD(v+\sqrt{\frac{D'}{D}})-b^2Du\sqrt{DD'}}
{a^2+abDu+v^{-1}b^2\sqrt{DD'}(1-Du^2)}
\notag
\end{align}
and 
\begin{align}\label{Eq:ValueofMCwhenclose}
\Phi_\varphi(t)&=\theta(a-b\sqrt{D})\psi\lrb{\varpi_\F^{-c(\theta)/e} \frac{ab}{a^2-b^2D}\lrb{\frac{D}{D'}v+v^{-1}(1-Du^2)+2\sqrt{\frac{D}{D'}}}}\\
&=\theta(a-b\sqrt{D})\psi\lrb{\varpi_\F^{-c(\theta)/e} \frac{ab}{(a^2-b^2D)v}\lrb{\lrb{\sqrt{\frac{D}{D'}}v+1}^2-Du^2}} .\notag
\end{align}
\end{lem}
\begin{proof}[Sketch of proof]
The congruence conditions guarantee that $k^{-1}tk\in \L^\times K_\mathfrak{A}(n)$, so \eqref{Eq:rewritingconjugation} is just to write $k^{-1}tk$ as a product of elements from $\L^\times$ and $K_\mathfrak{A}(n)$. For the value of matrix coefficient, we use Corollary \ref{Cor:MCofGeneralMinimalVec}, definition for $\tilde{\theta}$ in \eqref{Eq:thetatilde} and the special shape of $\alpha_\theta$ in \eqref{eq2.1:specialembedding}. In particular note that $$\zxz{a}{-b\sqrt{\frac{D}{D'}}}{-b\sqrt{DD'}}{a}=a-b\sqrt{\frac{D}{D'}} \cdot \sqrt{D'}=a-b\sqrt{D}$$ under the embedding of $\L$. 
\end{proof}

\begin{lem}\label{Lem:DomainofSupp}
Suppose that $l>0$ and $I(\varphi,\chi)\neq 0$. 
If $n-l$ is even, then $v_\F\lrb{\lrb{\sqrt{\frac{D}{D'}}v+1}^2}=n-l$, $v_\F(u)\geq (n-l)/2$(actually we can pick $u=0$) and $v_\F(b)\geq  \lfloor n/2\rfloor -\frac{n-l}{2}=\lfloor l/2\rfloor$. If $n-l$ is odd, then $v_\F\lrb{\lrb{\sqrt{\frac{D}{D'}}v+1}^2}>n-l$, $v_\F(u)=\frac{n-l-1}{2}$ and $v_\F(b)\geq \lceil n/2\rceil-\frac{n-l+1}{2}=\lceil \frac{l-1}{2}\rceil =\lfloor l/2\rfloor$.
\end{lem}
\begin{proof}
$I(\varphi,\chi)\neq 0$ in particular implies that $\Phi_\varphi(t)\chi(t)=1$ on  the support of the integral, on which we can apply the previous lemma. The level in $b$ of the part 
$$\psi\lrb{\varpi_\F^{-c(\theta)/e} \frac{ab}{(a^2-b^2D)v}\lrb{\lrb{\sqrt{\frac{D}{D'}}v+1}^2-Du^2}}$$
is $n-\min\{v_\F\lrb{\lrb{\sqrt{\frac{D}{D'}}v+1}^2}, v_\F(Du^2)\}$. One can then easily get the results on $v_\F\lrb{\sqrt{\frac{D}{D'}}v+1}$, $v_\F(u)$  by comparing this level with the level of $\overline{\theta}\chi(a+b\sqrt{D})$. The range for $b$ follows then from \eqref{Eq:Support1} and \eqref{Eq:Supp2}.

\end{proof}According to the range of $b$ in this lemma, we can further write 
$$\overline{\theta}\chi(a+b\sqrt{D})=\psi_\E\lrb{-\alpha_{\theta\overline{\chi}}\frac{b\sqrt{D}}{a}}=\psi\lrb{-2\alpha_{\theta\overline{\chi}}\frac{b\sqrt{D}}{a}}$$ where $v_\E(\alpha_{\theta\overline{\chi}})=-2l-1$. 
Then by Lemma \ref{lem:supportimpliesclosetorus} we can write 
\begin{align}
\Phi_\varphi(t)\chi(t)&= \psi\lrb{\varpi_\F^{-c(\theta)/e} \frac{ab}{(a^2-b^2D)v}\lrb{\lrb{\sqrt{\frac{D}{D'}}v+1}^2-Du^2}} \psi\lrb{-2\alpha_{\theta\overline{\chi}}\frac{b\sqrt{D}}{a}}\\
&=\psi\lrb{\varpi_\F^{-n} \frac{b}{av}\lrb{\lrb{\sqrt{\frac{D}{D'}}v+1}^2-Du^2}}\psi\lrb{-2\alpha_{\theta\overline{\chi}}\frac{b\sqrt{D}}{a}}\notag
\end{align}
The second equality follows from that $$\frac{ab}{(a^2-b^2D)v}=\frac{ab}{a^2v}(1+\frac{b^2D}{a^2}+\cdots).$$ One can show that all error terms do not matter by studying their valuations and using that $\psi$ is trivial on $\CO_\F$.

Then $\Phi_\varphi(t)\chi(t)=1$ is equivalent to that
\begin{equation}\label{eq:new-necessary-same-ramifiedAppendix}
b\lrb{\frac{D}{D'}v^2-\lrb{2\varpi_\F^n\alpha_{\theta\overline{\chi}}\sqrt{D}-2\sqrt{\frac{D}{D'}}}v+(1-Du^2)}\equiv 0\mod\varpi_\F^n
\end{equation}
on the support of the integral. 

The  discriminant of it is
\begin{equation}
\Delta(u)=4\varpi_\F^{n}\alpha_{\theta\overline{\chi}}\sqrt{D}\lrb{\varpi_\F^{n}\alpha_{\theta\overline{\chi}}\sqrt{D}-2\sqrt{\frac{D}{D'}}}+4\frac{D}{D'}Du^2.
\end{equation}

Consider first the case $n-l$ is even. We pick $u=0$ directly for simplicity. 
Then $v_\F(\Delta(0))=n-l$ and
\begin{equation}
{\Delta}(0)\equiv -8\varpi_\F^n\alpha_{\theta\overline{\chi}}\frac{D}{D'}\sqrt{D'}
\end{equation}
Recall that by Lemma \ref{lem:Epsilontestsameramified} (3), we write $\nu=\overline{\Theta}\chi\eta$, where $\eta $ is a character of $\E^\times$ extending $\eta_{\E/\F}$. In particular if we pick $\eta=\overline{\Delta_\theta}$, then $\nu=\overline{\theta}\chi$. When we write $\nu(1+\varpi_\E^{c(\nu)-1}x)=\theta(1+\varpi_\E^{c(\theta)-1}\delta x)$, this implies that
\begin{equation}
-\alpha_{\theta\overline{\chi}}\varpi_\E^{2l+1}\equiv \alpha_\theta \varpi_\E^{2n+1}\delta  \equiv    \varpi_\E\frac{1}{\sqrt{D'}} \delta \mod\varpi.
\end{equation}
Thus 
\begin{equation}
\Delta(0)\equiv 8\varpi_\F^{n-l} \frac{D}{D'} \delta .
\end{equation}
It's a square iff $2\delta $ is a square, consistent with Lemma \ref{lem:Epsilontestsameramified} (3).

When $\Delta(0)$ is indeed a square, we get two solutions of $v\mod \varpi_\F^{ \lceil n/2\rceil}$. 
For each of these two solutions we have by Lemma \ref{Lem:DomainofSupp}
\begin{equation}
I(\varphi,\chi)=\Vol(\BF^\times\backslash \BF^\times\{v_\F(a)=0, v_\F(b)\geq \lfloor l/2\rfloor\})=\frac{1}{q^{\lfloor l/2\rfloor}}.
\end{equation}

Now if $n-l$ is odd, $v_\F(\Delta(0))=n-l$ is odd, thus $\Delta(0)$ can never be a square. We need to pick $u$ such that $v_\F(u)=\frac{n-l-1}{2}$ and $\Delta(0)+4\frac{D}{D'}Du^2$ can be of higher evaluation. For this purpose we need that
\begin{equation}
8\varpi_\F^{n-l} \frac{D}{D'} \delta +4\frac{D}{D'}Du^2\equiv 0.
\end{equation}
Note that $D$ differs from $\varpi_\F$ by a square and $\frac{D}{D'}$ is also a square. So this being possible is equivalent to that $-2\delta $ is a square. This again is consistent with Lemma \ref{lem:Epsilontestsameramified} (3).

Once $-2\delta $ is a square, we can easily adjust $u$ so that $\Delta(u)$ is a square. In this case it's possible to get more solutions of $v\mod\varpi_\F^{\lceil n/2\rceil}$. For each solution we have by Lemma \ref{Lem:DomainofSupp}
\begin{equation}
I(\varphi,\chi)=\Vol(\BF^\times\backslash \BF^\times\{v_\F(a)=0, v_\F(b)\geq \lfloor l/2\rfloor\})=\frac{1}{q^{\lfloor l/2\rfloor}}.
\end{equation}


\subsubsection{Case $c(\theta\overline{\chi})=0$}
Note that $c(\theta\overline{\chi})=0$ implies that $\theta\overline{\chi}$ is an unramified character, not necessarily trivial as $\theta\overline{\chi}(\varpi_\E)=\pm 1$. 

Arguments in the previous case still apply, and we have on the $v_\F(b)\geq v_\F(a)=0$ part of the support of the integral
\begin{equation}
\Phi_\varphi(t)\chi(t)=\psi\lrb{\varpi_\F^{-c(\theta)/e} \frac{b}{av}\lrb{\lrb{\sqrt{\frac{D}{D'}}v+1}^2-Du^2}}.
\end{equation}

Thus 
\begin{equation}
b\lrb{\lrb{\sqrt{\frac{D}{D'}}v+1}^2-Du^2}\equiv 0\mod \varpi_\F^n.
\end{equation}

Recall from \eqref{Eq:Support1} and \eqref{Eq:Supp2} we have that $v_\F(b)\geq \max\{\lceil n/2\rceil -1-v_\F(u),\lfloor n/2\rfloor-v_\F(\sqrt{\frac{D}{D'}}v+1) \}$. Then we must have
$v_\F(u)\geq \lfloor n/2\rfloor$, $v_\F(\sqrt{\frac{D}{D'}}v+1)\geq \lceil n/2\rceil$. By the property of minimal vectors, we get essentially one test vector for such $u,v$. For simplicity we can pick $v=-\sqrt{\frac{D'}{D}}$ and $u=0$. Then the support of the integral is the whole torus. Now we just need to check $\Phi_\varphi(t)\chi(t)=1$ on $v_\F(a)>v_\F(b)=0$ part. But this part is just $\sqrt{D}$ times the part $ v_\F(b)\geq v_\F(a)=0$. So we just check 
\begin{align}
\Phi_{\varphi_0}\lrb{\zxz{-\sqrt{\frac{D}{D'}}}{0}{0}{1}\zxz{0}{1}{D}{0}\zxz{-\sqrt{\frac{D'}{D}}}{0}{0}{1}}\chi(\sqrt{D})&=\Phi_{\varphi_0}\lrb{-\zxz{0}{\sqrt{\frac{D}{D'}}}{\sqrt{DD'}}{0}}\chi(\sqrt{D})\\
&=\overline{\theta}\chi(\sqrt{D}).\notag
\end{align}
Since $\sqrt{D}$ differs from $\varpi_\E$ by an element from $\CO_\F^\times$, 
\begin{equation}
\overline{\theta}\chi(\sqrt{D})=1
\end{equation}
iff
\begin{equation}
\theta\overline{\chi}(\varpi_\E)=1.
\end{equation}
Thus we can find a test vector if and only if $\theta\overline{\chi}(\varpi_\E)=1$, and when this is true, for $u=0$ there exists a unique $v\mod\varpi_\F^{\lceil n/2\rceil}$  such that
\begin{equation}
I(\varphi,\chi)=2.
\end{equation}
Recall by our normalisation the volume of $\CO_\F^\times\backslash \CO_\E^\times$ is 1 and total volume of $\F^\times\backslash \E^\times$ is 2.

\bibliographystyle{alpha}
\bibliography{reference}

\begin{thebibliography}{CST17}

\bibitem[AL70]{AL1970}
A.O.L. Atkin and J.~Lehner.
\newblock Hecke operators on ${\Gamma}_0({N})$.
\newblock {\em Mathematische Annalen}, 185:134--160, 1970.

\bibitem[BH06]{BushnellHenniart:06a}
Colin~J. Bushnell and Guy Henniart.
\newblock {\em The local {L}anglands conjecture for {$\rm GL(2)$}}, volume 335
  of {\em Grundlehren der Mathematischen Wissenschaften [Fundamental Principles
  of Mathematical Sciences]}.
\newblock Springer-Verlag, Berlin, 2006.

\bibitem[Cas73]{Casselman1973}
William Casselman.
\newblock On some results of {A}tkin and {L}ehner.
\newblock {\em Mathematische Annalen}, 201:301--314, 1973.

\bibitem[CST14]{CST14}
Li~Cai, Jie Shu, and Ye~Tian.
\newblock Explicit {G}ross-{Z}agier and {W}aldspurger formulae.
\newblock {\em Algebra $\&$ Number Theory}, 8(10):2523--2572, 2014.

\bibitem[CST17]{CST17}
L.~Cai, J.~Shu, and Y.~Tian.
\newblock Cube sum problem and an explicit {G}ross-{Z}agier formula.
\newblock {\em Amer. Jour. of Math.}, 139(3):785--816, 2017.

\bibitem[DV09]{DV1}
Samit Dasgupta and John Voight.
\newblock Heegner points and sylvester's conjecture.
\newblock {\em Arithmetic Geometry: Clay Mathematics Institute Summer School,
  Arithmetic Geometry, July 17-August 11, 2006, Mathematisches Institut,
  Georg-August-Universit{\"a}t, G{\"o}ttingen, Germany}, 8:91, 2009.

\bibitem[DV17]{DV17}
S.~Dasgupta and J.~Voight.
\newblock Sylvester's problem and mock {H}eegner points.
\newblock {\em arXiv:1707.05874}, 2017.

\bibitem[Gro88]{Gross88}
Benedict~H. Gross.
\newblock Local orders, root numbers, and modular curves.
\newblock {\em American Journal of Mathematics}, 110(6):1153--1182, 1988.

\bibitem[GZ86]{GZ1986}
B.H. Gross and D.B. Zagier.
\newblock Heegner points and derivatives of {L}-series.
\newblock {\em Inventiones mathematicae}, 84:225--320, 1986.

\bibitem[HN]{HN18}
Yueke Hu and Paul Nelson.
\newblock New test vector for waldspurger's period integral.
\newblock {\em Preprint}.

\bibitem[KM88]{KM1988}
M.~A. Kenku and Fumiyuki Momose.
\newblock Automorphism groups of the modular curves ${X}_0({N})$.
\newblock {\em Compositio Mathematica}, 65(1):51--80, 1988.

\bibitem[Kob13]{Koba2013}
S.~Kobayashi.
\newblock The $p$-adic {G}ross-{Z}agier formula for elliptic curves at
  supersingular primes.
\newblock {\em Invent. Math.}, 191:527--629, 2013.

\bibitem[Kol90]{Kolyvagin1990}
V.~A. Kolyvagin.
\newblock Euler systems.
\newblock In {\em The Grothendieck Festschrift II}, volume~87 of {\em Progr.
  Math.}, pages 435--483. Birkhauser Boston, Boston, MA, 1990.

\bibitem[Lan]{Langlands}
R.P. Langlands.
\newblock On the functional equation of the artin l-functions.
\newblock {\em https://publications.ias.edu/sites/default/files/a-ps.pdf}.

\bibitem[Liv95]{Liverance}
E.~Liverance.
\newblock A formula for the root number of a family of elliptic curves.
\newblock {\em Journal of Number Theory}, 51(2):288 -- 305, 1995.

\bibitem[LLT]{LLT}
Yongxiong Li, Yu~Liu, and Ye~Tian.
\newblock On the {B}irch and {S}winnerton-{D}yer conjecture for {CM} elliptic
  curves over $\mathbb{Q}$.
\newblock {\em arXiv:1605.01481}.

\bibitem[Neu99]{Neukirchbook1}
J\"{u}rgen Neukirch.
\newblock {\em Algebraic number theory}, volume 322 of {\em Grundlehren der
  Mathematischen Wissenschaften [Fundamental Principles of Mathematical
  Sciences]}.
\newblock Springer-Verlag, Berlin, 1999.
\newblock Translated from the 1992 German original and with a note by Norbert
  Schappacher, With a foreword by G. Harder.

\bibitem[Ogg80]{Ogg80}
A.~P. Ogg.
\newblock Modular functions.
\newblock In {\em Santa Cruz Conference on Finite Groups}, volume~37 of {\em
  Proc. Sympos. Pure Math.}, pages 521--532. Amer. Math. Soc., Providence,
  1980.

\bibitem[PR87]{PR1987}
B.~Perrin-Riou.
\newblock Points de heegner et deriv{\'e}es de fonctions {L} p-adiques.
\newblock {\em Inventiones mathematicae}, 89:455--510, 1987.

\bibitem[Sai93]{Saito:93a}
Hiroshi Saito.
\newblock On {T}unnell's formula for characters of ${GL}(2)$.
\newblock {\em Compositio Mathematica}, 85(1):99--108, 1993.

\bibitem[Sat86]{Satge}
P.~Satg{\'e}.
\newblock Groupes de {S}elmer et corps cubiques.
\newblock {\em J. Number Theory}, 23(3):294--317, 1986.

\bibitem[Sel51]{Selmer51}
E.~S. Selmer.
\newblock The diophatine equation $ax^3+by^3+cz^3=0$.
\newblock {\em Acta Math.}, 87:203--362, 1951.

\bibitem[Shi94]{Shimurabook}
Goro Shimura.
\newblock {\em Introduction to the arithmetic theory of automorphic functions},
  volume~11 of {\em Publications of the Mathematical Society of Japan}.
\newblock Princeton University Press, Princeton, NJ, 1994.
\newblock Reprint of the 1971 original, Kan{\^o} Memorial Lectures, 1.

\bibitem[Sil94]{Silvermanbook2}
Joseph~H. Silverman.
\newblock {\em Advanced topics in the arithmetic of elliptic curves}, volume
  151 of {\em Graduate Texts in Mathematics}.
\newblock Springer-Verlag, New York, 1994.

\bibitem[Syl79]{Sylv}
J.~J. Sylvester.
\newblock On certain tenary cubic-form equations.
\newblock {\em Amer. J. Math.}, 2(4):357--393, 1879.

\bibitem[Tun83]{Tunnell:83a}
Jerrold~B. Tunnell.
\newblock Local {$\epsilon $}-factors and characters of {${\rm GL}(2)$}.
\newblock {\em Amer. J. Math.}, 105(6):1277--1307, 1983.

\bibitem[YZZ13]{YZZ}
Xinyi Yuan, Shou-Wu Zhang, and Wei Zhang.
\newblock {\em The {G}ross-{Z}agier formula on {S}himura curves}, volume 184 of
  {\em Annals of Mathematics Studies}.
\newblock Princeton University Press, Princeton, NJ, 2013.

\bibitem[ZK87]{ZK}
D.~Zagier and G.~Kramarz.
\newblock Numerical investigations related to the {$L$}-series of certain
  elliptic curves.
\newblock {\em J. Indian Math. Soc. (N.S.)}, 52:51--69 (1988), 1987.

\end{thebibliography}
\end{document}